\documentclass[12pt]{article}
\usepackage{mathrsfs}
\usepackage{epic,eepic,epsf,epsfig}
\usepackage{amsfonts,srcltx,mathrsfs}
\usepackage{multirow}
\usepackage{amsmath}
\usepackage{amssymb}
\usepackage{amsbsy}
\usepackage{graphicx}
\usepackage{amsfonts}
\usepackage{color}
\usepackage{setspace}
\newtheorem{lem}{Lemma}[section]%
\newtheorem{theorem}[lem]{Theorem}%
\newtheorem{cor}[lem]{Corollary}%

\def\nd{\mathrel{\bigm|\kern-.7em/}}

\def\f{\noindent}

\def\P\GammaL{\hbox{\rm P\GammaL}}

\def\mod{\hbox{\rm mod }}

\begin{document}
\title{Spectral skeletons and applications (an updated version)}

\footnotetext{E-mails: zhangwq@pku.edu.cn}

\author{Wenqian Zhang\\
{\small School of Mathematics and Statistics, Shandong University of Technology}\\
{\small Zibo, Shandong 255000, P.R. China}}
\date{}
\maketitle

\begin{abstract}
For a graph $G$, its spectral radius $\rho(G)$ is the largest eigenvalue of its adjacency matrix. Let $\mathcal{F}$ be a finite family of graphs with $\min_{F\in \mathcal{F}}\chi(F)=r+1\geq3$, where $\chi(F)$ is the chromatic number of $F$.  Set $t=\max_{F\in\mathcal{F}}|F|$. Let $T(rt,r)$ be the Tur\'{a}n graph of order $rt$ with $r$ parts. Assume that some $F_{0}\subseteq\mathcal{F}$ is a subgraph of the graph obtained from $T(rt,r)$ by embedding a path or a matching in one part. Let ${\rm EX}(n,\mathcal{F})$ be the set of graphs with the maximum number of edges among all the graphs of order $n$ containing not any $F\in\mathcal{F}$. Simonovits \cite{S1,S2} gave general results on the graphs in ${\rm EX}(n,\mathcal{F})$. Let ${\rm SPEX}(n,\mathcal{F})$ be the set of graphs with the maximum spectral radius among all the graphs of order $n$ containing not any $F\in\mathcal{F}$. Motivated by the work of Simonovits, we characterize the specified structure of the graphs in ${\rm SPEX}(n,\mathcal{F})$ in this paper. Moreover, some applications are also included.

\bigskip

\f {\bf Keywords:} spectral radius; spectral extremal graph; extremal graph; symmetric subgraph.

\end{abstract}

\baselineskip 17 pt

\section{Introduction}

All graphs considered in this paper are finite, undirected and without multi-edges or loops. We first introduce some elementary notations in graph theory. For a graph $G$, let $\overline{G}$ denote its complement.  The vertex set and edge set of $G$ are denoted by $V(G)$ and $E(G)$, respectively. Let $|G|=|V(G)|$ and $e(G)=|E(G)|$. For a subset $S\subseteq V(G)$, let $G[S]$ be the subgraph of $G$ induced by $S$, and let $G-S=G[V(G)-S]$. Let $E(S)$ be the set of edges of $G$ inside $S$ and let $e(S)=|E(S)|$. For a vertex $u$,  another vertex $v$ is called a {\em neighbor} of $u$, if they are adjacent in $G$.  let $N_{G}(u)$ be the set of neighbors of $u$ in $G$, and let $N_{S}(u)$ be the set of neighbors of $u$ in $S$. The degree of $u$ is $d_{G}(u)=|N_{G}(u)|$. Let $\delta(G)$ denote the minimum degree of $G$. For any terminology used but not defined here, one may refer to \cite{CRS}.

For $\ell\geq2$  graphs $G_{1},G_{2},...,G_{\ell}$, let $\cup_{1\leq i\leq \ell}G_{i}$ (or $G_{1}\cup G_{2}\cup\cdots\cup G_{\ell}$) be the disjoint union of them. Let $\prod_{1\leq i\leq \ell}G_{i}$ (or $G_{1}\prod G_{2}\prod\cdots\prod G_{\ell}$) be the graph obtained from $\cup_{1\leq i\leq \ell}G_{i}$ by connecting each vertex in $G_{i}$ to each vertex in $G_{j}$ for any $1\leq i\neq j\leq \ell$. For a positive integer $t$ and a graph $H$, let $tH=\cup_{1\leq i\leq \ell}H$. For a certain integer $n$, let $K_{n}, C_{n}$ and $P_{n}$ be the complete graph, the cycle and the path of order $n$, respectively. For $r\geq2$, let $K_{n_{1},n_{2},...,n_{r}}$ be the complete $r$-partite graph with parts of sizes $n_{1},n_{2},...,n_{r}$. Let $T(n,r)$ be the Tur\'{a}n graph of order $n$ with $r$ parts (i.e., the complete $r$-partite graph of order $n$, in which each part has $\lfloor\frac{n}{r}\rfloor$ or $\lceil\frac{n}{r}\rceil$ vertices). Assume that $T(n,1)=\overline{K_{n}}$, and $T(n,0)$ is an empty graph.

 Let $G$ be a graph with vertices $v_{1},v_{2},...,v_{n}$. Its {\em adjacency matrix} is defined as $A(G)=(a_{ij})_{n\times n}$, where $a_{ij}=1$ if  $v_{i}$ and $v_{j}$ are adjacent, and $a_{ij}=0$ otherwise. The {\em spectral radius} $\rho(G)$ of $G$ is the largest eigenvalue of $A(G)$. By the Perron--Frobenius Theorem,  $\rho(G)$ has a non-negative eigenvector (called Perron vector). Moreover, $G$ has a positive eigenvector if $G$ is connected. For a graph $F$, denote by $F\subseteq G$ if $G$ contains a copy of $F$, and $F\nsubseteq G$ otherwise. For a set of graphs $\mathcal{F}$, $G$ is call $\mathcal{F}$-free, if $F\nsubseteq G$  for any $F\in \mathcal{F}$.   Let ${\rm EX}(n,\mathcal{F})$ be the set of $\mathcal{F}$-free graphs of order $n$ with the maximum number of edges, and let ${\rm SPEX}(n,\mathcal{F})$ be the set of $\mathcal{F}$-free graphs of order $n$ with the maximum spectral radius. Let ${\rm ex}(n,\mathcal{F})$ be the number of edges of any graph in ${\rm EX}(n,\mathcal{F})$.

The classical Tur\'{a}n Theorem \cite{T} states that ${\rm EX}(n,\left\{K_{r+1}\right\})=\left\{T(n,r)\right\}$. For a graph $G$, let $\chi(G)$ denote its {\em chromatic number}. Erd\H{o}s,
Stone and Simonovits \cite{E Simonovits,E Stone} proposed the stability theorem
$${\rm ex}(n,\mathcal{F})=(1-\frac{1}{\chi(\mathcal{F})-1})\frac{n^{2}}{2}+o(n^{2}),$$
 where $\chi(\mathcal{F})=\min_{F\in\mathcal{F}}\chi(F)$.
There are a lot of articles on Tur\'{a}n
type problem (for example, see \cite{B,CGPW,DJ,EFGG,HQL,HLF,L,YZ,YZ,Y,Y1,Y2,Y3}). In 2010, Nikiforov \cite{N1} presented formally a spectral version of Tur\'{a}n-type problem: what is the maximum spectral radius of graphs with specified subgraph structure?
 In recent years, this problem has attracted much attention from graph theorists (for example, see \cite{BDT,CFTZ,CDT1,CDT2,CDT3,DKLNTW,FLSZ,FTZ,LL,LP1,LP2,LZZ,NWK,WKX,WNKF,ZHL,ZL,Z}).

Simonovits \cite{S1} introduced the concept of symmetric subgraphs. Let $G$ be a graph. For $\tau\geq1$ induced subgraphs of $G$: $Q_{1},Q_{2},..,Q_{\tau}$, they are called symmetric subgraphs in $G$, if they are connected and vertex-disjoint, and there are isomorphisms $\psi_{j}: Q_{1}\rightarrow Q_{j}$ for all $2\leq j\leq \tau$ such that for any $u\in V(Q_{1})$ and $v\in V(G)-(\cup_{1\leq i\leq\tau}V(Q_{i}))$, $uv\in E(G)$ if and only if $\psi_{j}(u)v\in E(G)$. These symmetric subgraphs are called trivial if $|Q_{i}|=1$, and non-trivial otherwise. Trivial symmetric subgraphs are also called symmetric vertices. Clearly, $\tau$ vertices $v_{1},v_{2},...,v_{\tau}$ of $G$ are symmetric if and only if there are no edges inside $\left\{v_{1},v_{2},...,v_{\tau}\right\}$ and $N_{G}(v_{1})=N_{G}(v_{i})$ for any $1\leq i\leq\tau$.

\medskip

\f{\bf Observation.} Let $F$ and $G$ be two graphs, where $t=|F|\geq2$ and $F\nsubseteq G$. Assume that $Q_{1},Q_{2},...,Q_{t}$ are symmetric subgraphs of $G$. Let $G'$ be the graph obtained from $G$ by adding a new copy of $Q_{1}$, say $Q$, such that $Q,Q_{1},Q_{2},...,Q_{t}$ are symmetric subgraphs in $G'$. Then $F\nsubseteq G'$. In fact, if $F\subseteq G'$, then $V(F)\cap V(Q)\neq\emptyset$ as $F\nsubseteq G$. Noting $t=|F|$, there is at least one of $Q_{1},Q_{2},...,Q_{t}$, say $Q_{t}$, such that $V(F)\cap V(Q_{t})=\emptyset$. Let $F'$ be the subgraph of $G$ induced by $(V(F)-V(Q))\cup V(Q_{t})$. By symmetry of $Q$ and $Q_{t}$ in $G'$, we see that $F\subseteq F'$, which contradicts the fact that $F\nsubseteq G$. Hence $F\nsubseteq G'$ is obtained.

\medskip

The above observation is one of the most direct motivations to study symmetric subgraphs in extremal graph theory. For certain integers $n,r$ and $c$,
let $\mathbb{D}(n,r,c)$ denote the family of graphs $G$ of order $n$ satisfying the
following symmetry condition:\\
$(i)$  It is possible to omit at most $c$ vertices of $G$ so that the remaining graph $G'$ is of form $\prod_{1\leq i\leq r}G_{i}$, where $||V(G_{i})-\frac{n}{r}|\leq c$ for any $1\leq i\leq r$.\\
$(ii)$ For any $1\leq i\leq r$, $G_{i}=\cup_{1\leq j\leq k_{i}}H^{i}_{j}$, where $H^{i}_{1},H^{i}_{2},...,H^{i}_{k_{i}}$ are symmetric subgraphs in $G$.

Simonovits \cite{S1,S2} gave the following two results on extremal graphs.

\begin{theorem}{\rm (\cite{S1})}\label{Simonovits 1}
Let $\mathcal{F}$ be a finite family of graphs with $\min_{F\in \mathcal{F}}\chi(F)=r+1\geq3$. Set $t=\max_{F\in\mathcal{F}}|F|$. Assume that $F\subseteq P_{t}\prod T(t(r-1),r-1)$ for some $F\in\mathcal{F}$. Then, for sufficiently large $n$, $\mathbb{D}(n,r,c)$ contains a graph $G$ in ${\rm EX}(n,\mathcal{F})$, where $c$ is constant with respect to $n$. Furthermore, if $G$ is the only extremal graph in $\mathbb{D}(n,r,c)$, then it is
unique in ${\rm EX}(n,\mathcal{F})$.
\end{theorem}

\begin{theorem}{\rm (\cite{S2})}\label{Simonovits 2}
Let $\mathcal{F}$ be a finite family of graphs with $\min_{F\in \mathcal{F}}\chi(F)=r+1\geq3$. Set $t=\max_{F\in\mathcal{F}}|F|$. Assume that $F\subseteq tK_{2}\prod T(t(r-1),r-1)$ for some $F\in\mathcal{F}$. Then, for sufficiently large $n$, $\mathbb{D}(n,r,c)$ contains a graph in ${\rm EX}(n,\mathcal{F})$, where $c$ is constant with respect to $n$. Furthermore, any graph $G\in\mathbb{D}(n,r,c)\cap {\rm EX}(n,\mathcal{F})$ has the following property:\\
$(i)$ $G'=T(|G'|,r)$, where $G'$ is the graph obtained from $G$ by omitting at most $c$ vertices;\\
$(ii)$ each exceptional vertex in $V(G)-V(G')$ is connecting either to all the vertices of $G'$ or to all
the vertices of $r-1$ classes of $G'$ and to no vertex of the remaining class.
\end{theorem}

We should notice that there are probably other extremal graphs which are not in $\mathbb{D}(n,r,c)$ in Theorem \ref{Simonovits 1} and Theorem \ref{Simonovits 2}. In this paper, we will give two related spectral skeletons (see Theorem \ref{Pl-free} and Theorem \ref{matching-free}). Interestingly, all the spectral extremal graphs must be of the forms stated in the theorems. We first give the following result on the minimum degree of the spectral extremal graphs.

\begin{lem}\label{mini degree}
Let $\mathcal{F}$ be a finite family of graphs with $\min_{F\in \mathcal{F}}\chi(F)=r+1\geq3$. For every $\theta>0$, there exists $n_{0}$ such that if $G\in {\rm SPEX}(n,\mathcal{F})$ with $n\geq n_{0}$, then $G$ is connected and $\delta(G)>(\frac{r-1}{r}-\theta)n$.
\end{lem}

Now we state the two spectral skeletons as follows.

\begin{theorem}\label{Pl-free}
Let $\mathcal{F}$ be a finite family of graphs with $\min_{F\in \mathcal{F}}\chi(F)=r+1\geq3$. Set $t=\max_{F\in\mathcal{F}}|F|$. Assume that $F_{1}\subseteq P_{t}\prod T(t(r-1),r-1)$ and $F_{2}\subseteq tK_{2}\prod tK_{2}\prod T(t(r-2),r-2)$ for some $F_{1},F_{2}\in\mathcal{F}$. Then there exists a constant $M_{0}$ such that, for sufficiently large $n$, any graph $G$ in ${\rm SPEX}(n,\mathcal{F})$ has a partition $V(G)=\cup_{1\leq i\leq r}S_{i}$ satisfying:\\
 $(i)$ for any $1\leq i\leq r$, $|S_{i}|\rightarrow\infty$ when  $n\rightarrow\infty$;\\
 $(ii)$ for any $2\leq i\leq r$, there exists $S'_{i}\subseteq S_{i}$ such that $|S_{i}-S'_{i}|\leq M_{0}$, and $N_{G}(v)=V(G)-S_{i}$ for any $v\in S'_{i}$.
\end{theorem}

For a finite family of graphs $\mathcal{F}$ with $\min_{F\in \mathcal{F}}\chi(F)=r+1\geq3$ and $t=\max_{F\in\mathcal{F}}|F|$, define $q(\mathcal{F})$ to be the minimum integer $s$ such that $F\subseteq \overline{K_{s}}\prod T(tr,r)$ for some $F\subseteq\mathcal{F}$. Clearly, $q(\mathcal{F})\leq t$.

\begin{theorem}\label{matching-free}
Let $\mathcal{F}$ be a finite family of graphs with $\min_{F\in \mathcal{F}}\chi(F)=r+1\geq3$. Set $t=\max_{F\in\mathcal{F}}|F|$ and $q=q(\mathcal{F})$. Assume that $F_{0}\subseteq tK_{2}\prod T(t(r-1),r-1)$ for some $F_{0}\in\mathcal{F}$. Then there exists a constant $M_{0}$ such that, for sufficiently large $n$, any graph $G$ in ${\rm SPEX}(n,\mathcal{F})$ has a partition $V(G)=W\cup(\cup_{1\leq i\leq r}S_{i})$ satisfying the following property:\\
$(i)$ $|W|=q-1$ and $|S_{i}|=\lfloor\frac{n-q+1}{r}\rfloor$ or $\lceil\frac{n-q+1}{r}\rceil$ for any $1\leq i\leq r$;\\
$(ii)$ for any $1\leq i\leq r$, there exists $S'_{i}\subseteq S_{i}$ such that $|S_{i}-S'_{i}|\leq M_{0}$, and $N_{G}(v)=V(G)-S_{i}$ for any $v\in S'_{i}$.
\end{theorem}

Clearly, all the spectral extremal graphs in Theorem \ref{matching-free} are contained in $\mathbb{D}(n,r,c)$, where $c=O(1)$.
Assume that $1\leq s\leq\frac{m}{2}$ and $r\geq2$. As an old result, Simonovits (see \cite{S} or \cite{S1}) showed that, for sufficiently large $n$, $K_{s-1}\prod T(n-s+1,r)$ is the unique graph in ${\rm EX}(n,\left\{F\right\})$, where $F=(sK_{2}\cup\overline{K_{m-2s}})\prod T(m(r-1),r-1)$.
As an application of Theorem \ref{matching-free}, we give the spectral version of this result.

\begin{theorem}\label{application 1}
For integers $1\leq s\leq\frac{m}{2}$ and $r\geq2$, let $F=(sK_{2}\cup\overline{K_{m-2s}})\prod T(m(r-1),r-1)$. For large $n$, $K_{s-1}\prod T(n-s+1,r)$ is the unique graph in ${\rm SPEX}(n,\left\{F\right\})$.
\end{theorem}

Let $G$ be a graph. For two vertices $u$ and $v$, let $d_{G}(u,v)$ be the {\em distance} of them in $G$ (i.e., the shortest length of a path from $u$ to $v$). For an integer $k\geq2$, the {\em $k$-th power} of $G$, denoted by $G^{k}$, is defined as $V(G^{k})=V(G)$, and any two distinct vertices $u,v$ are adjacent in $G^{k}$ if and only if $d_{G}(u,v)\leq k$. Very recently, Fang and Zhao \cite{FZ} characterized the graphs in ${\rm EX}(n,\left\{C^{2}_{m}\right\})$ and ${\rm SPEX}(n,\left\{C^{2}_{m}\right\})$, where $m\geq6$ and $m\not\equiv0~(\mod3)$. (The case $m\equiv0~(\mod3)$ is left therein.) In this paper, we extends this result on $C^{2}_{m}$ to $C^{k}_{m}$, where $m\not\equiv0~(\mod(k+1))$ and $k\geq2$.

\begin{theorem}\label{application 2}
 Assume that $m=p(k+1)+h$, where $k,p\geq2$ and $1\leq h\leq k$ are integers. Let $s=h-\lfloor\frac{h-1}{p}\rfloor p$ and $r=k+\lceil\frac{h}{p}\rceil$.  For sufficiently large $n$, we have the following conclusions.\\
 $(i)$ $K_{s-1}\prod T(n-s+1,r)$ is the unique graph in ${\rm EX}(n,\left\{C^{k}_{m}\right\})$.\\
 $(ii)$ $K_{s-1}\prod T(n-s+1,r)$ is the unique graph in ${\rm SPEX}(n,\left\{C^{k}_{m}\right\})$.
\end{theorem}

For the (spectral) extremal graphs on $C^{k}_{n}$ in the case of $m=p(k+1)$, we make the following conjecture.

\medskip

\f{\bf Conjecture.}  Assume that $m=p(k+1)$, where $k\geq2$ and $p\geq3$ are integers. For sufficiently large $n$, each graph $G$ in ${\rm EX}(n,\left\{C^{k}_{m}\right\})$ or ${\rm SPEX}(n,\left\{C^{k}_{m}\right\})$ is of form $H\prod T(n-|H|,k-1)$, where $H$ is a $\left\{C_{2p}\right\}$-free graph.

\medskip

The rest of this paper is organized as follows. In Section 2, we include some lemmas needed in the proofs of the main results. In Section 3, we give the proof of Lemma \ref{mini degree}. The proofs of Theorem \ref{Pl-free} and Theorem \ref{matching-free} are given in Section 4. In Section 5, we give the proofs of Theorem \ref{application 1} and Theorem \ref{application 2}.

\section{Preliminaries}

 To prove the main results of this paper, we first include several lemmas.

\begin{lem}{\rm (\cite{CRS})}\label{subgraph}
If $H$ is a subgraph of a connected graph $G$, then $\rho(H)\leq\rho(G)$, with equality if and only if $H=G$.
\end{lem}

The following lemma is a variation (with a very similar proof) of Theorem 8.1.3 in \cite{CRS}.

\begin{lem}{\rm (\cite{CRS})}\label{eigenvector trans}
Let $G$ be a connected graph with a Perron vector $\mathbf{x}=(x_{w})_{w\in V(G)}$. For a vertex $u\in V(G)$, Let $G'$ be the graph obtained from $G$ by deleting edges $uv_{1},uv_{2},...,uv_{s}$, and adding edges $uw_{1},uw_{2},...,uw_{t}$, where $s,t\geq1$. If $\sum_{1\leq j\leq t}x_{w_{j}}\geq\sum_{1\leq i\leq s}x_{v_{i}}$ and $\left\{v_{1},v_{2},...,v_{s}\right\}\neq\left\{w_{1},w_{2},...,w_{s}\right\}$, then $\rho(G')>\rho(G)$.
\end{lem}

The following lemma will be used to show that the vertex partition in Theorem \ref{matching-free} is balanced. It may have its own interests.

\begin{lem} \label{eigenvector trans var}
Let $G$ be a connected graph with a Perron vector $\mathbf{x}=(x_{w})_{w\in V(G)}$, where $\sum_{w\in V(G)}x^{2}_{w}=1$. Assume that the vertex set $V(G)$ can be partitioned into three parts $U,A,B$, such that $|U|>0,|A|>|B|>0$, and the following holds.\\
$(i)$ $G[A\cup B]$ is a complete bipartite graph with parts $A$ and $B$;\\
$(ii)$ there are two subsets $U_{1},U_{2}\subseteq U$, such that $N_{U}(v)=U_{1}$ for any $v\in A$, and $N_{U}(w)=U_{2}$ for any $w\in B$;\\
 $(iii)$ $\sum_{v\in U_{1}}x_{v}<\sum_{v\in U_{2}}x_{v}$.\\
Let $G'$ be the graph obtained from $G$ by deleting all the edges between $A$ and $U_{1}$, and between $B$ and $U_{2}$, and adding all the edges between $A$ and $U_{2}$, and between $B$ and $U_{1}$. Then $\lambda_{1}(G')>\lambda_{1}(G)$.
\end{lem}

\f{\bf Proof:} Let $A(G)$ and $A(G')$ be the adjacency matrices of $G$ and $G'$, respectively. Let $\mathbf{x}^{T}$ be the transpose of $\mathbf{x}$. Recall that $\sum_{v\in V(G)}x^{2}_{v}=1$. By symmetry,  we can assume that $x_{v}=a>0$ for any $v\in A$, and $x_{w}=b>0$ for any $w\in B$. Set $\sum_{v\in U_{1}}x_{v}=u_{1}$ and $\sum_{v\in U_{2}}x_{v}=u_{2}$. Then $u_{1}<u_{2}$. Set $n_{1}=|A|$ and $n_{2}=|B|$. Then $n_{1}>n_{2}\geq 1$.

\f{\bf Case 1.} $an_{1}>bn_{2}$.

Since
\begin{equation}
\begin{aligned}
&\mathbf{x}^{T}A(G')\mathbf{x}-\mathbf{x}^{T}A(G)\mathbf{x}\\
&=2(\sum_{uv\in E(G')}x_{u}x_{v})-2(\sum_{uv\in E(G)}x_{u}x_{v})\\
&=2(an_{1}u_{2}+bn_{2}u_{1})-2(an_{1}u_{1}+bn_{2}u_{2})\\
&=2(an_{1}-bn_{2})(u_{2}-u_{1})\\
&>0,
\end{aligned}\notag
\end{equation}
we see that $\mathbf{x}^{T}A(G')\mathbf{x}>\mathbf{x}^{T}A(G)\mathbf{x}=\rho(G)$. Then $\rho(G')>\rho(G)$ by linear algebra.

\f{\bf Case 2.} $an_{1}\leq bn_{2}$.

 Recall that $n_{1}>n_{2}$. Then $b\sqrt{n_{2}}\geq\frac{an_{1}}{\sqrt{n_{2}}}>a\sqrt{n_{1}}$. Since $u_{1}<u_{2}$, it follows that $u_{2}b\sqrt{n_{2}}-u_{1}a\sqrt{n_{1}}>0$.
Now define a new vector $\mathbf{y}=(y_{v})_{v\in V(G)}$, where
\begin{equation}
y_{v}=\left\{
\begin{array}{ll}
b\sqrt{\frac{n_{2}}{n_{1}}}& for~v\in A,\\
a\sqrt{\frac{n_{1}}{n_{2}}}& for~v\in B,\\
x_{v}& for~v\in V(G)-(A\cup B).
\end{array}
\right.\notag
\end{equation}
Clearly, $\sum_{v\in V(G)}y^{2}_{v}=\sum_{v\in V(G)}x^{2}_{v}=1$.
Since
\begin{equation}
\begin{aligned}
&\mathbf{y}^{T}A(G')\mathbf{y}-\mathbf{x}^{T}A(G)\mathbf{x}\\
&=2(\sum_{uv\in E(G')}y_{u}y_{v})-2(\sum_{uv\in E(G)}x_{u}x_{v})\\
&=2\left((u_{2}n_{1}b\sqrt{\frac{n_{2}}{n_{1}}}+u_{1}n_{2}a\sqrt{\frac{n_{1}}{n_{2}}}+
n_{1}b\sqrt{\frac{n_{2}}{n_{1}}}\cdot n_{2}a\sqrt{\frac{n_{1}}{n_{2}}}\right)-2(u_{1}an_{1}+u_{2}bn_{2}+an_{1}bn_{2})\\
&=2u_{2}b(\sqrt{n_{1}n_{2}}-n_{2})-2u_{1}a(n_{1}-\sqrt{n_{1}n_{2}})\\
&=2(\sqrt{n_{1}}-\sqrt{n_{2}})(u_{2}b\sqrt{n_{2}}-u_{1}a\sqrt{n_{1}})\\
&>0,
\end{aligned}\notag
\end{equation}
we see that $\mathbf{y}^{T}A(G')\mathbf{y}>\mathbf{x}^{T}A(G)\mathbf{x}=\rho(G)$. Then $\rho(G')>\rho(G)$ by linear algebra.
The proof is completed. \hfill$\Box$

\medskip

The following result is the spectral version of the Stability Lemma due to Nikiforov \cite{N2}.

\begin{lem} {\rm (\cite{N2})}\label{spec-stability}
Suppose $r\geq2,\frac{1}{\ln n}<c<r^{-8(r+1)(r+21)}$ and $0<\epsilon<2^{-36}r^{-24}$. Let $G$ be a graph of order $n$. If $\rho(G)>(\frac{r-1}{r}-\epsilon)n$, then one of the following holds:\\
$(i)$ $G$ contains a complete $(r+1)$-partite graph $K_{\lfloor c\ln n\rfloor,\lfloor c\ln n\rfloor,...,\lfloor c\ln n\rfloor,\lceil n^{1-\sqrt{c}}\rceil}$;\\
$(ii)$ $G$  differs from $T(n,r)$ in fewer than $(\epsilon^{\frac{1}{4}}+c^{\frac{1}{8r+8}})n^{2}$ edges.
\end{lem}

The following result is a direct corollary of Lemma \ref{spec-stability} (similar to \cite{DKLNTW}).

\begin{cor} \label{var-spec-stability}
Let $\mathcal{F}$ be a finite family of graphs with $\min_{F\in \mathcal{F}}\chi(F)=r+1\geq3$. For every $\epsilon>0$, there exist $\delta>0$ and $n_{0}$ such that if $G$ is an $\mathcal{F}$-free graph of order $n\geq n_{0}$ with $\rho(G)>(\frac{r-1}{r}-\delta)n$, then $G$ can be
obtained from $T(n,r)$ by adding and deleting at most $\epsilon n^{2}$ edges.
\end{cor}

The following lemma was given in \cite{CFTZ}.

\begin{lem}\label{sets inequality}
Let $A_{1},A_{2},...,A_{\ell}$ be $\ell\geq2$ finite sets. Then
$$|\cap_{1\leq i\leq \ell}A_{i}|\geq(\sum_{1\leq i\leq\ell}|A_{i}|)-(\ell-1)|\cup_{1\leq i\leq \ell}A_{i}|.$$
\end{lem}

The following result is taken from \cite{EG}

\begin{lem}\label{ex n Pt}
Let $n\geq k\geq2$. Then
${\rm ex}(n,\left\{P_{k}\right\})\leq\frac{k-2}{2}n$.
\end{lem}

The following result is an easy consequence of Lemma 3.2.1 in \cite{S1} (just by choosing $m_{k}=\tau$ for any $k\geq1$ therein). Furthermore, Simonovits pointed out that the orders of the symmetric subgraphs are bounded (i.e., $O(1)$) at the end of the proof of Lemma 3.2.1 in \cite{S1}.

\begin{cor}\label{symmetric}
Let $G$ be a $\left\{P_{\ell}\right\}$-free graph of order $n$. For any integer $\tau>0$, there exist two positive integers $C(\ell,\tau)$ and $D(\ell,\tau)$ such that if $n\geq C(\ell,\tau)$, then $G$ has $\tau$ symmetric subgraphs $Q_{1},Q_{2},...,Q_{\tau}$, where $|Q_{i}|\leq D(\ell,\tau)$ for any $1\leq i\leq \tau$.
\end{cor}

\section{Proof of Lemma \ref{mini degree}}

Let $\mathcal{F}$ be a finite family of graphs with $\min_{F\in \mathcal{F}}\chi(F)=r+1\geq3$. Set $t=\max_{F\in\mathcal{F}}|F|$.
Let $G\in {\rm SPEX}(n,\mathcal{F})$, where $n$ is sufficiently large. Let $\eta>0$ be a small constant with respect to $n$, and $\eta\leq\frac{\theta}{4}$. By the result of Erd\H{o}s, Stone and Simonovits \cite{E Simonovits,E Stone}, we have the following fact.

\medskip

\f{\bf Fact 1.} ${\rm ex}(n',\mathcal{F})\leq(\frac{r-1}{2r}+\eta^{2})n'^{2}$ for any large $n'$.

\medskip

\begin{lem}\label{Lemma spec-lower}
 $\rho(G)\geq\frac{r-1}{r}n-\frac{r}{4n}$.
 \end{lem}

\f{\bf Proof:} Note that $T(n,r)$ is $\mathcal{F}$-free, since $\min_{F\in \mathcal{F}}\chi(F)=r+1$. Let $\overline{d}$ denote the average degree of $T(n,r)$. As is well known, $\rho(T(n,r))\geq\overline{d}\geq\frac{r-1}{r}n-\frac{r}{4n}$. Noting that $G\in {\rm SPEX}(n,\mathcal{F})$, we have
$\rho(G)\geq\rho(T(n,r))\geq\frac{r-1}{r}n-\frac{r}{4n},$ as desired. \hfill$\Box$

\begin{lem}\label{Lem connected}
 $G$ is connected.
 \end{lem}

\f{\bf Proof:} Suppose that $G$ is not connected. Let $G_{1}$ be a component of $G$ with $\rho(G_{1})=\rho(G)$. Let $n_{1}=|G_{1}|<n$, and let $\Delta$ be the maximum degree of $G_{1}$. Then $\Delta\geq\rho(G_{1})=\rho(G)\geq(\frac{r-1}{r}-\eta)n$ by Lemma \ref{Lemma spec-lower}. Thus $\Delta>t$ as $n$ is large. Let $u$ be a vertex of $G_{1}$ with degree $\Delta$. Assume that $v$ is a vertex in $V(G)-V(G_{1})$. Let $G'$ be the graph obtained from $G$ by deleting the edges incident with $v$ and adding a new edge $uv$. By Lemma \ref{subgraph}, we have $\rho(G')>\rho(G_{1})=\rho(G)$. This implies that $G'$ is not $\mathcal{F}$-free. Then $F\subseteq G'$ for some $F\in\mathcal{F}$. Since $u$ has degree $\Delta>t$ in $G$, there is a neighbor $u_{1}$ of $u$ such that $u_{1}\notin V(F)$. Let $F'$ be the subgraph of $G$ induced by $(V(F)-\left\{v\right\})\cup\left\{u_{1}\right\}$. Clearly, $F\subseteq F'$. This contradicts the fact that $G$ is $\mathcal{F}$-free. Hence $G$ is connected. \hfill$\Box$

\begin{lem}\label{Lem partition}
$e(G)\geq(\frac{r-1}{2r}-\eta^{3})n^{2}$. Moreover, there is a partition $V(G)=\cup_{1\leq i\leq r}V_{i}$ such that $\sum_{1\leq i\leq r}e(G[V_{i}])$ is minimum, where $\sum_{1\leq i\leq r}e(G[V_{i}])\leq\eta^{3}n^{2}$ and $||V_{i}|-\frac{n}{r}|\leq\eta n$ for any $1\leq i\leq r$.
\end{lem}

\f{\bf Proof:} Since $G$ is $\mathcal{F}$-free and $\rho(G)\geq\frac{r-1}{r}n-\frac{r}{4n}$, by Corollary \ref{var-spec-stability} (letting $\epsilon=\frac{1}{2}\eta^{3}$), we see that $G$ can be obtained from $T(n,r)$ by deleting and adding at most $\frac{1}{2}\eta^{3}n^{2}$ edges. It follows that $e(G)\geq(\frac{r-1}{2r}-\eta^{3})n^{2}$, and there is a (balanced) partition $V(G)=\cup_{1\leq i\leq r}U_{i}$ such that $\sum_{1\leq i\leq r}e(G[U_{i}])\leq\eta^{3}n^{2}$. Now we select a partition $V(G)=\cup_{1\leq i\leq r}V_{i}$ such that $\sum_{1\leq i\leq r}e(G[V_{i}])$ is minimum. Then
 $$\sum_{1\leq i\leq r}e(G[V_{i}])\leq\sum_{1\leq i\leq r}e(G[U_{i}])\leq\eta^{3}n^{2}.$$

Let $a=\max_{1\leq i\leq r}||V_{i}|-\frac{n}{r}|$. Without loss of generality, assume that $a=||V_{1}|-\frac{n}{r}|$. Using the Cauchy-Schwarz inequality, we see
$$2\sum_{2\leq i<j\leq r}|V_{i}||V_{j}|=(\sum_{2\leq i\leq r}|V_{i}|)^{2}-\sum_{2\leq i\leq r}|V_{i}|^{2}\leq\frac{r-2}{r-1}(\sum_{2\leq i\leq r}|V_{i}|)^{2}=\frac{r-2}{r-1}(n-|V_{1}|)^{2}.$$
Then, noting $a=||V_{1}|-\frac{n}{r}|$,
\begin{equation}
\begin{aligned}
e(G)&\leq(\sum_{1\leq i<j\leq r}|V_{i}||V_{j}|)+(\sum_{1\leq i\leq r}e(G[V_{i}]))\\
&\leq|V_{1}|(n-|V_{1}|)+(\sum_{2\leq i<j\leq r}|V_{i}||V_{j}|)+\eta^{3}n^{2}\\
&\leq|V_{1}|(n-|V_{1}|)+\frac{r-2}{2(r-1)}(n-|V_{1}|)^{2}+\eta^{3}n^{2}\\
&=\frac{r-1}{2r}n^{2}-\frac{r}{2(r-1)}a^{2}+\eta^{3}n^{2}.
\end{aligned}\notag
\end{equation}
Recall that $e(G)\geq(\frac{r-1}{2r}-\eta^{3})n^{2}$. Thus, $(\frac{r-1}{2r}-\eta^{3})n^{2}\leq\frac{r-1}{2r}n^{2}-\frac{r}{2(r-1)}a^{2}+\eta^{3}n^{2}$, implying that $a\leq\sqrt{\frac{4(r-1)}{r}\eta^{3}n^{2}}\leq\eta n$. This finishes the proof. \hfill$\Box$

\medskip

Let $L=\left\{u\in V(G)~|~d_{G}(u)\leq(\frac{r-1}{r}-4\eta)n\right\}$.

\begin{lem}\label{Lem L bound}
$|L|\leq\eta n$.
\end{lem}
\f{\bf Proof:} Suppose $|L|>\eta n$. Let $L'\subseteq L$ with $|L'|=\lfloor\eta n\rfloor$. Set $n'=|G-L'|=n-\lfloor\eta n\rfloor< (1-\eta )n+1$. Recall $e(G)\geq(\frac{r-1}{2r}-\eta^{3})n^{2}$. Then
\begin{equation}
\begin{aligned}
e(G-L')&\geq e(G)-\sum_{w\in L'}d_{G}(w)\\
&\geq(\frac{r-1}{2r}-\eta^{3})n^{2}-(\frac{r-1}{r}-4\eta)n\eta n\\
&=(\frac{r-1}{2r}-\frac{r-1}{r}\eta+4\eta^{2}-\eta^{3})n^{2}\\
&=\frac{r-1}{2r}((1-\eta)n)^{2}+(4\eta^{2}-\eta^{3}-\frac{r-1}{2r}\eta^{2})n^{2}\\
&>\frac{r-1}{2r}((1-\eta)n+1)^{2}+\eta^{2}n^{2}\\
&>(\frac{r-1}{2r}+\eta^{2})((1-\eta)n+1)^{2}\\
&>(\frac{r-1}{2r}+\eta^{2})n'^{2}.
\end{aligned}\notag
\end{equation}
Since $G-L'$ is $\mathcal{F}$-free, by Fact 1, we have that ${\rm ex}(n',\mathcal{F})\leq(\frac{r-1}{2r}+\eta^{2})n'^{2}$ as $n'$ is large. But this contradicts that $e(G-L')>(\frac{r-1}{2r}+\eta^{2})n'^{2}$.  Hence $|L|\leq\eta n$. \hfill$\Box$

\medskip

Let $W_{i}=\left\{v\in V_{i}~|~d_{V_{i}}(v)\geq2\eta n\right\}$ for $1\leq i\leq r$, and let $W=\cup_{1\leq i\leq r}W_{i}$.

\begin{lem}\label{Lem W bound}
$|W|\leq\eta^{2} n$.
\end{lem}
\f{\bf Proof:} Since $\sum_{1\leq i\leq r}e(G[V_{i}])\leq\eta^{3}n^{2}$, and
$$\sum_{1\leq i\leq r}e(G[V_{i}])=\sum_{1\leq i\leq r}\frac{1}{2}\sum_{v\in V_{i}}d_{V_{i}}(v)\geq\sum_{1\leq i\leq r}\frac{1}{2}\sum_{v\in W_{i}}d_{V_{i}}(v)\geq\sum_{1\leq i\leq r}\eta n|W_{i}|=\eta n|W|,$$
we have $\eta n|W|\leq\eta^{3}n^{2}$, implying that $|W|\leq\eta^{2}n$.
 \hfill$\Box$

\medskip

For any $1\leq i\leq r$, let $\overline{V}_{i}=V_{i}-(L\cup W)$.

\begin{lem}\label{Lem common neighbors}
Let $1\leq \ell\leq r$ be fixed. Assume that $u_{1},u_{2},...,u_{t}\in \cup_{1\leq i\neq \ell\leq r}\overline{V}_{i}$. Then there are $t+1$ vertices in $V_{\ell}$ which are the common neighbors of $u_{1},u_{2},...,u_{t}$ in $G$.
\end{lem}

\f{\bf Proof:} For any $1\leq i\leq t$, assume that $u_{i}\in \overline{V}_{j_{i}}$, where $j_{i}\neq\ell$. Then $d_{G}(u_{i})\geq(\frac{r-1}{r}-4\eta)n$ and $d_{V_{j_{i}}}(u_{i})\leq2\eta n$ as $u_{i}\notin L\cup W$. Recall that $|V_{s}|\leq\frac{n}{r}+\eta n$ for any $1\leq s\leq r$.
Hence $$d_{V_{\ell}}(u_{i})= d_{G}(u_{i})-d_{V_{j_{i}}}(u_{i})-\sum_{1\leq s\neq\ell,j_{i}\leq r}d_{V_{s}}(u_{i})\geq(\frac{1}{r}-(r+4)\eta)n.$$
Then, by Lemma \ref{sets inequality}, we have
$$|\cap_{1\leq i\leq t}N_{V_{\ell}}(u_{i})|\geq(\sum_{1\leq i\leq t}|N_{V_{\ell}}(u_{i})|)-(t-1)|V_{\ell}|\geq(\frac{1}{r}-(r+5)t\eta)n\geq t+1.$$
Thus, there are $t+1$ vertices in $V_{\ell}$ which are adjacent to all the vertices $u_{1},u_{2},...,u_{t}$.
 \hfill$\Box$

\begin{lem} \label{Lem L empty}
$L=\emptyset$.
\end{lem}
\f{\bf Proof:} Let $\mathbf{x}=(x_{v})_{v\in V(G)}$ be a Perron vector of $G$ such that $\max_{v\in V(G)}x_{v}=1$. Assume that $u^{*}$ is a vertex with $x_{u^{*}}=1$. Let $v^{*}$ be a vertex in $V(G)-W$ such that $x_{v^{*}}=\max_{v\in V(G)-W}x_{v}$. Recall that $|W|\leq\eta^{2}n$ by Lemma \ref{Lem W bound}.
Since $\rho(G)=\rho(G)x_{u^{*}}\leq|W|+(n-|W|)x_{v^{*}}$, we have $$x_{v^{*}}\geq\frac{\rho(G)-|W|}{n-|W|}\geq\frac{\rho(G)-|W|}{n}\geq(\frac{r-1}{r}-2\eta^{2})>\frac{1}{3}.$$
Then  $|W|\leq3\eta^{2}nx_{v^{*}}$.

Assume that $v^{*}\in V_{i_{0}}$. Then $|N_{\overline{V}_{i_{0}}}(v^{*})|\leq|N_{V_{i_{0}}}(v^{*})|\leq2\eta n$ as $v^{*}\notin W$.
Hence
\begin{equation}
\begin{aligned}
\rho(G)x_{v^{*}}&=(\sum_{v\in N_{L\cup W}(v^{*})}x_{v})+(\sum_{v\in N_{\overline{V}_{i_{0}}}(v^{*})}x_{v})
+(\sum_{v\in N_{\cup_{1\leq i\neq i_{0}\leq r}\overline{V}_{i}}(v^{*})}x_{v})\\
&\leq(|W|+|L|x_{v^{*}})+2\eta nx_{v^{*}}+\sum_{v\in \cup_{1\leq i\neq i_{0}\leq r}\overline{V}_{i}}x_{v}\\
&\leq(3\eta+3\eta^{2})nx_{v^{*}}+\sum_{v\in \cup_{1\leq i\neq i_{0}\leq r}\overline{V}_{i}}x_{v}.
\end{aligned}\notag
\end{equation}
It follows that
$$\sum_{v\in \cup_{1\leq i\neq i_{0}\leq r}\overline{V}_{i}}x_{v}\geq(\rho(G)-3\eta n-3\eta^{2}n)x_{v^{*}}.$$

Suppose that $L\neq\emptyset$. Let $u\in L$. Then $d_{G}(u)\leq(\frac{r-1}{r}-4\eta)n<\rho(G)-4\eta n+\eta^{2}n$.
Thus $$\sum_{v\in N_{G}(u)}x_{v}=(\sum_{v\in N_{W}(u)}x_{v})+(\sum_{v\in N_{V(G)-W}(u)}x_{v})\leq |W|+d_{G}(u)x_{v^{*}}\leq(\rho(G)-4\eta n+4\eta^{2}n)x_{v^{*}}.$$
It follows that
 $$(\sum_{v\in \cup_{1\leq i\neq i_{0}\leq r}\overline{V}_{i}}x_{v})-(\sum_{v\in N_{G}(u)}x_{v})\geq(\eta n-7\eta^{2}n)x_{v^{*}}>0.$$

Let $G'$ be the graph obtained
from $G$ by deleting all the edges incident with $u$ and adding all the edges between $u$ and $\cup_{1\leq i\neq i_{0}\leq r}\overline{V}_{i}$.
Then
$$\mathbf{x}^{T}(\rho(G')-\rho(G))\mathbf{x}\geq\mathbf{x}^{T}(A(G')-A(G))\mathbf{x}=2x_{u}((\sum_{v\in \cup_{1\leq i\neq i_{0}\leq r}\overline{V}_{i}}x_{v})-(\sum_{v\in N_{G}(u)}x_{v}))>0,$$
which implies that $\rho(G')>\rho(G)$.
Since $G\in {\rm SPEX}(n,\mathcal{F})$, we have that $G'$ is not $\mathcal{F}$-free. Let $F\in\mathcal{F}$ be a subgraph of $G'$. If $u\notin V(F)$, then $F\subseteq G$, a contradiction. So, we can assume that $u\in V(F)$. Let $u_{1},u_{2},...,u_{s}$ be all the neighbors of $u$ in $F$, where $s\leq t$. Note that $u_{1},u_{2},...,u_{s}\in\cup_{1\leq i\neq i_{0}\leq r}\overline{V}_{i}$.
By Lemma \ref{Lem common neighbors}, there are $t+1$ vertices in $V_{i_{0}}$, which are the common neighbors of $u_{1},u_{2},...,u_{s}$ in $G$. Hence, there is a vertex $w\notin V(F)$, which is adjacent to all the vertices  $u_{1},u_{2},...,u_{s}$ in $G$. Let $F'$ be the subgraph of $G$ induced by $(V(F)-\left\{u\right\})\cup\left\{w\right\}$. Clearly, $F\subseteq F'$. This contradicts the fact that $G$ is $\mathcal{F}$-free. Hence $L=\emptyset$. This completes the proof.
 \hfill$\Box$

\medskip

\f{\bf Proof of Lemma \ref{mini degree}.} Clearly, Lemma \ref{mini degree} follows from Lemma \ref{Lem connected} and Lemma \ref{Lem L empty}. \hfill$\Box$

\section{Proofs of Theorem \ref{Pl-free} and Theorem \ref{matching-free}}

Let $\mathcal{F}$ be a finite family of graphs with $\min_{F\in \mathcal{F}}\chi(F)=r+1\geq3$.  Set $t=\max_{F\in\mathcal{F}}|F|$. Assume that $\mathcal{F}$ satisfies property $(P)$: there are some graphs $F_{1},F_{2}\in\mathcal{F}$ such that
$$F_{1}\subseteq P_{t}\prod T(t(r-1),r-1)$$
 and
$$F_{2}\subseteq tK_{2}\prod tK_{2}\prod T(t(r-2),r-2).$$
Note that if $F\subseteq tK_{2}\prod T(t(r-1),r-1)$ for some $F\in\mathcal{F}$, then $\mathcal{F}$ satisfies property $(P)$. So, property $(P)$ is satisfied in both Theorem \ref{Pl-free} and Theorem \ref{matching-free}.

For sufficiently large $n$, let $G\in {\rm SPEX}(n,\mathcal{F})$, where $\mathcal{F}$ satisfies the above property $(P)$. Let $\theta>0$ be a small constant with respect to $n$ (such that all the following inequalities on it are satisfied). By Lemma \ref{mini degree}, we have that $G$ is connected and $\delta(G)\geq(\frac{r-1}{r}-\theta)n$.  Similar to Lemma \ref{Lemma spec-lower}, the following bound holds:
$$\rho(G)\geq\frac{r-1}{r}n-\frac{r}{4n}.$$
Since $G$ is $\mathcal{F}$-free, by Corollary \ref{var-spec-stability} (letting $\epsilon=\frac{1}{2}\theta^{3}$), we see that $G$ can be obtained from $T(n,r)$ by deleting and adding at most $\frac{1}{2}\theta^{3}n^{2}$ edges. It follows that
 $$(\frac{r-1}{2r}-\theta^{3})n^{2}\leq e(G)\leq(\frac{r-1}{2r}+\theta^{3})n^{2}.$$
  Let $V(G)=\cup_{1\leq i\leq r}V_{i}$ be a partition such that $\sum_{1\leq i\leq r}e(G[V_{i}])$ is minimum. Similar to Lemma \ref{Lem partition}, we have that $\sum_{1\leq i\leq r}e(G[V_{i}])\leq\theta^{3}n^{2}$ and $||V_{i}|-\frac{n}{r}|\leq\theta n$ for any
  $1\leq i\leq r$.

For any $1\leq i\leq r$, let $W_{i}=\left\{v\in V_{i}~|~d_{V_{i}}(v)\geq\theta n\right\}$, and let $\overline{V}_{i}=V_{i}-W_{i}$.
Set $W=\cup_{1\leq i\leq r}W_{i}$.

\begin{lem}\label{Main W bound}
 $|W|\leq M$, where $M=\lceil4^{(r-1)t}(\frac{2}{\theta})^{t}rt\rceil$ is constant with respect to $n$.
\end{lem}
\f{\bf Proof:} Similar to Lemma \ref{Lem W bound}, we can show that $|W|\leq2\theta^{2}n$.
For any fixed $1\leq i\leq r$ and $1\leq j\neq i\leq r$, we have $d_{V_{i}}(v)\leq d_{V_{j}}(v)$ for any $v\in V_{i}$, since $\sum_{1\leq \ell\leq r}e(G[V_{i}])$ is minimum (otherwise, we can obtain a contradiction by moving $v$ from $V_{i}$ to $V_{j}$). Therefore, $d_{V_{j}}(v)\geq\frac{1}{2}(d_{V_{i}}(v)+d_{V_{j}}(v))$ for any $v\in V_{i}$.  Recall that $|V_{s}|\leq(\frac{1}{r}+\theta)n$ for any $1\leq s\leq r$.
If $v\in W_{i}$, then
 $$d_{V_{j}}(v)\geq \frac{1}{2}(d_{G}(v)-\sum_{1\leq s\neq i,j\leq r}|V_{s}|)\geq \frac{1}{2}(\frac{r-1}{r}n-\theta n-(r-2)(\frac{1}{r}+\theta)n)\geq(\frac{1}{2r}-\frac{r-1}{2}\theta)n.$$
 It follows that
 $$d_{\overline{V}_{j}}(v)\geq d_{V_{j}}(v)-|W|\geq(\frac{1}{2r}-\frac{r-1}{2}\theta)n-2\theta^{2}n\geq\frac{n}{2r}-\frac{r}{2}\theta n.$$
If $v\in \overline{V}_{i}$, then $d_{V_{i}}(v)\leq \theta n$.
Thus
$$d_{V_{j}}(v)\geq d_{G}(v)-d_{V_{i}}(v)-\sum_{1\leq s\neq i,j\leq r}|V_{s}|\geq(\frac{r-1}{r}-\theta)n-\theta n-(r-2)(\frac{1}{r}+\theta)n\geq(\frac{1}{r}-r\theta)n.$$
It follows that
 $$d_{\overline{V}_{j}}(v)\geq d_{V_{j}}(v)-|W|\geq(\frac{1}{r}-r\theta)n-2\theta^{2}n\geq\frac{n}{r}-(r+1)\theta n.$$

In the following, we first show that $|W_{1}|\leq4^{(r-1)t}(\frac{2}{\theta})^{t}t$.

For any $v\in W_{1}$, we have that $d_{\overline{V}_{1}}(v)\geq d_{V_{1}}(v)-|W_{1}|\geq\theta n-2\theta^{2}n\geq\frac{1}{2}\theta n$. For $j\neq1$, we have $d_{\overline{V}_{j}}(v)\geq\frac{n}{2r}-\frac{r}{2}\theta n$ by the above discussion.
 Let $$Y_{1}=\left\{(w,Z)~|~w\in W_{1},Z\subseteq\overline{V}_{1},|Z|=t\right\},$$
  where $w$ is adjacent to the all the vertices in $Z$.
Since $d_{\overline{V}_{1}}(w)\geq\frac{1}{2}\theta n$ for any $w\in W_{1}$,  we have $|Y_{1}|\geq\tbinom{\frac{1}{2}\theta n}{t}|W_{1}|$. On the other hand,  $\overline{V}_{1}$ has at most $\tbinom{(\frac{1}{r}+\theta) n}{t}$ subsets $Z$ of size $t$. Hence, there is a subset $Z_{1}\subseteq \overline{V}_{1}$ of size $t$ such that the vertices in $Z_{1}$ have at least $\frac{\tbinom{\frac{1}{2}\theta n}{t}|W_{1}|}{\tbinom{(\frac{1}{r}+\theta) n}{t}}\geq(\frac{\theta}{2})^{t}|W_{1}|$ common neighbors  in $W_{1}$. Let $W^{1}\subseteq W_{1}$ be the set of common neighbors of the vertices in $Z_{1}$. Then $|W^{1}|\geq(\frac{\theta}{2})^{t}|W_{1}|$.

Suppose that for some $1\leq s\leq r-1$, we have obtained several sets $Z_{i}\subseteq \overline{V}_{i}$ of size $t$ for any $1\leq i\leq s$ and $W^{s}\subseteq W_{1}$ , such that any two vertices from distinct parts of $W^{s},Z_{1},...,Z_{s}$ are adjacent in $G$ and $|W^{s}|\geq(\frac{1}{4})^{(s-1)t}(\frac{\theta}{2})^{t}|W_{1}|$. Let $V'_{s+1}\subseteq \overline{V}_{s+1}$ be the set of common neighbors of the vertices in $\cup_{1\leq i\leq s}Z_{i}$. Since $(\frac{n}{r}-2\theta)n\leq|\overline{V}_{s+1}|\leq(\frac{n}{r}+2\theta)n$ and $d_{\overline{V}_{s+1}}(u)\geq\frac{n}{r}-(r+1)\theta n$ for any $u\in\cup_{1\leq i\leq s}Z_{i}$ by the above discussion,
we have that $d_{\overline{V}_{s+1}}(u)\geq|\overline{V}_{s+1}|-(r+3)\theta n$ for any $u\in\cup_{1\leq i\leq s}Z_{i}$. By Lemma \ref{sets inequality}, we have
$$|V'_{s+1}|=|\bigcap_{u\in\cup_{1\leq i\leq s}Z_{i}}N_{\overline{V}_{s+1}}(u)|\geq|\overline{V}_{s+1}|-st(r+3)\theta n\geq|\overline{V}_{s+1}|-rt(r+3)\theta n.$$
For any $w\in W_{1}$, note that $d_{\overline{V}_{s+1}}(w)\geq\frac{n}{2r}-\frac{r}{2}\theta n$ as above. Then $d_{V'_{s+1}}(w)\geq\frac{n}{2r}-rt(r+4)\theta n$.

Note that $|V'_{s+1}|\leq(\frac{1}{r}+\theta)n$.
Let $$Y_{s+1}=\left\{(w,Z)~|~w\in W^{s},Z\subseteq V'_{s+1},|Z|=t\right\},$$
 where $w$ is adjacent to all the vertices in $Z$.
Since $d_{V'_{s+1}}(w)\geq\frac{n}{2r}-rt(r+4)\theta n\geq\frac{n}{3r}$ for any $w\in W^{s}$,  we have $|Y_{s+1}|\geq|W^{s}|\tbinom{\frac{n}{3r}}{t}$. On the other hand,  $V'_{s+1}$ has at most $\tbinom{(\frac{1}{r}+\theta) n}{t}$ subsets $Z$ of size $t$. Hence, there is a subset $Z_{s+1}\subseteq V'_{s+1}$ of size $t$ such that the vertices in $Z_{s+1}$ have at least $\frac{|W^{s}|\tbinom{\frac{n}{3r}}{t}}{\tbinom{(\frac{1}{r}+\theta) n}{t}}\geq(\frac{1}{4})^{t}|W^{s}|$ common neighbors  in $W^{s}$. Let $W^{s+1}\subseteq W^{s}$ be the set of common neighbors of the vertices in $Z_{s+1}$. Then $$|W^{s+1}|\geq(\frac{1}{4})^{t}|W^{s}|\geq(\frac{1}{4})^{st}(\frac{\theta}{2})^{t}|W_{1}|.$$
Therefore, by induction on $s$, we can obtain several sets $Z_{i}\subseteq \overline{V_{i}}$ of size $t$ for any $1\leq i\leq r$  and $W^{r}\subseteq W_{1}$, such that any two vertices from distinct parts of $W^{r},Z_{1},...,Z_{r}$ are adjacent in $G$ and $|W^{r}|\geq(\frac{1}{4})^{(r-1)t}(\frac{\theta}{2})^{t}|W_{1}|$. So, if $|W_{1}|\geq4^{(r-1)t}(\frac{2}{\theta})^{t}t$, then $|W^{r}|\geq t$. It follows that the subgraph induced by $W^{r},Z_{1},...,Z_{r}$ contains a complete $(r+1)$-partite graph of which each part has size at least $t$. This contradicts the fact that $G$ is $\mathcal{F}$-free. Thus $|W_{1}|<4^{(r-1)t}(\frac{2}{\theta})^{t}t$. Similarly,  we can show that $|W_{i}|<4^{(r-1)t}(\frac{2}{\theta})^{t}t$ for any $2\leq i\leq r$. Hence $|W|<4^{(r-1)t}(\frac{2}{\theta})^{t}rt$.
This finishes the proof. \hfill$\Box$

\medskip

Similar to Lemma \ref{Lem common neighbors}, the following lemma also holds.

\begin{lem}\label{Main neighbors}
Let $1\leq \ell\leq r$ be fixed. Assume that $u_{1},u_{2},...,u_{2rt}\in \cup_{1\leq i\neq \ell\leq r}\overline{V}_{i}$. Then there are $t+1$ vertices in $\overline{V}_{\ell}$ which are adjacent to all the vertices $u_{1},u_{2},...,u_{2rt}$ in $G$.
\end{lem}

Let $\mathbf{x}=(x_{v})_{v\in V(G)}$ be a Perron vector of $G$ such that $\max_{v\in V(G)}x_{v}=1$. Assume that $u^{*}$ is a vertex with $x_{u^{*}}=1$. Let $v^{*}$ be a vertex in $V(G)-W$ such that $x_{v^{*}}=\max_{v\in V(G)-W}x_{v}$.

\begin{lem}\label{entry lower}
For any $u\in V(G)$, we have $x_{u}\geq (1-5\theta)x_{v^{*}}$ and $x_{u}\geq\frac{r-1}{r}-6\theta$.
\end{lem}

\f{\bf Proof:} Recall that $|W|\leq M$.
Since $\rho(G)=\rho(G)x_{u^{*}}\leq|W|+(n-|W|)x_{v^{*}}$, we have $$x_{v^{*}}\geq\frac{\rho(G)-|W|}{n}\geq\frac{r-1}{r}-\theta.$$
It follows that $x_{v^{*}}>\frac{1}{3}$ and $|W|\leq3Mx_{v^{*}}$.

Assume that $v^{*}\in V_{i_{0}}$, where $1\leq i_{0}\leq r$. Then $|N_{\overline{V}_{i_{0}}}(v^{*})|\leq|N_{V_{i_{0}}}(v^{*})|\leq\theta n$ as $v^{*}\notin W$.
Hence
\begin{equation}
\begin{aligned}
\rho(G)x_{v^{*}}&=(\sum_{v\in N_{ W}(v^{*})}x_{v})+(\sum_{v\in N_{\overline{V}_{i_{0}}}(v^{*})}x_{v})+(\sum_{v\in N_{\cup_{1\leq i\neq i_{0}\leq r}\overline{V}_{i}}(v^{*})}x_{v})\\
&\leq|W|+\theta nx_{v^{*}}+\sum_{v\in \cup_{1\leq i\neq i_{0}\leq r}\overline{V}_{i}}x_{v}\\
&\leq(3M+\theta n)x_{v^{*}}+\sum_{v\in \cup_{1\leq i\neq i_{0}\leq r}\overline{V}_{i}}x_{v}.
\end{aligned}\notag
\end{equation}
Noting that $\rho(G)\geq\frac{n}{3}$, it follows that
$$\sum_{v\in \cup_{1\leq i\neq i_{0}\leq r}\overline{V}_{i}}x_{v}\geq(\rho(G)-3M-\theta n)x_{v^{*}}\geq(1-4\theta)\rho(G)x_{v^{*}}.$$

Suppose that $x_{u}< (1-5\theta)x_{v^{*}}$ for some $u\in V(G)$.  Then
$$x_{u}+\sum_{v\in N_{G}(u)}x_{v}=(1+\rho(G))x_{u}<(1-4\theta)\rho(G)x_{v^{*}}$$
So, $$\sum_{v\in \cup_{1\leq i\neq i_{0}\leq r}\overline{V}_{i}}x_{v}>x_{u}+\sum_{v\in N_{G}(u)}x_{v}.$$
Let $G'$ be the graph obtained
from $G$ by deleting all the edges incident with $u$ and adding all the edges between $u$ and $\cup_{1\leq i\neq i_{0}\leq r}\overline{V}_{i}$.
Then
$$\mathbf{x}^{T}(\rho(G')-\rho(G))\mathbf{x}\geq\mathbf{x}^{T}(A(G')-A(G))\mathbf{x}=2x_{u}((\sum_{v\in \cup_{1\leq i\neq i_{0}\leq r}\overline{V}_{i}}x_{v})-(x_{u}+\sum_{v\in N_{G}(u)}x_{v}))>0,$$
implying that $\rho(G')>\rho(G)$.
Then $F\subseteq G'$ for some $F\in \mathcal{F}$, since $G\in {\rm SPEX}(n,\mathcal{F})$. Similar to Lemma \ref{Lem L empty}, we will obtain a contradiction. Hence we have $x_{u}\geq (1-5\theta)x_{v^{*}}$ for any $u\in V(G)$. Since $x_{v^{*}}\geq\frac{r-1}{r}-\theta$ as above, we see that $x_{u}\geq\frac{r-1}{r}-6\theta$.
This finishes the proof.
 \hfill$\Box$

\medskip

For any fixed $1\leq i\leq r$, let $Q$ be an induced subgraph of $G[\overline{V}_{i}]$ (or a subset of $\overline{V}_{i}$).
Define
$$\alpha(Q)=E(Q)\cup\left\{uv\in E(G)~|~u\in V(Q),v\in (V_{i}-V(Q))\cup W\right\},$$
and $$\beta(Q)=\left\{uv\notin E(G)~|~u\in V(Q),v\in\cup_{1\leq j\neq i\leq r}\overline{V}_{j}\right\}.$$

\begin{lem}\label{symmetry edge bound}
For any fixed $1\leq i_{0}\leq r$, if $Q_{1},Q_{2},...,Q_{t}$ are symmetric subgraphs of $G[\overline{V}_{i_{0}}]$, then $|\alpha(Q_{i})|\leq(M+2t)|Q_{1}|$ and $|\beta(Q_{1})|\leq9(M+2t)|Q_{1}|$.
\end{lem}

\f{\bf Proof:} For convenience, assume $i_{0}=1$. We first show that there are at most $t$ vertices in $\overline{V}_{1}-\cup_{1\leq i\leq t}V(Q_{i})$ which have neighbors in $V(Q_{1})$. If not, there are $t$ vertices $u_{1},u_{2},...,u_{t}$ in $\overline{V}_{1}-\cup_{1\leq i\leq t}V(Q_{i})$, such that $u_{i}v_{i}$ is an edge for any $1\leq i\leq t$, where $v_{1},v_{2},...,v_{t}\in V(Q_{1})$ (may be repeated). Since $Q_{1},Q_{2},...,Q_{t}$ are symmetric subgraphs in $G[\overline{V}_{1}]$, for $2\leq j\leq t$ let $v^{j}_{1},v^{j}_{2},...,v^{j}_{t}\in V(Q_{j})$ be the symmetric vertices of $v_{1},v_{2},...,v_{t}$ in $V(Q_{j})$, respectively. Then $u_{i}$ is adjacent to $v^{j}_{i}$ for any $1 \leq i,j\leq t$. Note that $Q_{j}$ is connected for any $1\leq j\leq t$. Let $P_{v^{j}_{j}v^{j}_{j+1}}$ be a path connecting $v^{j}_{j}$ to $v^{j}_{j+1}$ in $Q_{j}$. Then
$$u_{1}P_{v^{1}_{1}v^{1}_{2}}u_{2}P_{v^{2}_{2}v^{2}_{3}}u_{3}\cdots u_{t-1}P_{v^{t-1}_{t-1}v^{t-1}_{t}}u_{t}$$
is a path of order at least $2t-1$ in $G[\overline{V}_{1}]$. Now we select its a sub-path of order $t$, say $w_{1},w_{2},...,w_{t}$. By Lemma \ref{Main neighbors}, these $t$ vertices have $t$ common neighbors in $\overline{V}_{2}$, say $w^{2}_{1},w^{2}_{2},...,w^{2}_{t}$. By Lemma \ref{Main neighbors}, these $2t$ vertices have $t$ common neighbors in $\overline{V}_{3}$. Repeat this process. Eventually, we can obtain a subgraph in $G$, which contains $P_{t}\prod T(t(r-1),r-1)$ as a subgraph, a contradiction. Hence there are at most $t$ vertices in $\overline{V}_{1}-\cup_{1\leq i\leq t}V(Q_{i})$ which have neighbors in $V(Q_{1})$. Moreover, using a similar discussion, we see that $Q_{i}$ is $\left\{P_{t}\right\}$-free for any $1\leq i\leq t$. This implies that $e(Q_{1})\leq t|Q_{1}|$ by Lemma \ref{ex n Pt}. Hence, $|\alpha(Q_{1})|\leq (M+2t)|Q_{1}|$.

Suppose $|\beta(Q_{1})|>9(M+2t)|Q_{1}|$. Let $G'$ be the graph obtained from $G$ by deleting all the edges in $\alpha(Q_{1})$, and adding all the non-edges in $\beta(Q_{1})$.
By Lemma \ref{entry lower}, we see that $x_{u}\geq\frac{r-1}{r}-6\theta\geq\frac{1}{3}$ for any $u\in V(G)$. Then
$$(\sum_{uv\in \beta(Q_{1})}x_{u}x_{v})-(\sum_{uv\in \alpha(Q_{1})}x_{u}x_{v})\geq(\frac{1}{3})^{2}|\beta(Q_{1})|-|\alpha(Q_{1})|>0.$$
It follows that
$$\mathbf{x}^{T}(\rho(G')-\rho(G))\mathbf{x}\geq\mathbf{x}^{T}(A(G')-A(G))\mathbf{x}=2(\sum_{uv\in \beta(Q_{1})}x_{u}x_{v})-2(\sum_{uv\in \alpha(Q_{1})}x_{u}x_{v})>0,$$
implying that $\rho(G')>\rho(G)$.
Then $F\subseteq G'$ for some $F\in\mathcal{F}$ as $G\in {\rm SPEX}(n,\mathcal{F})$.
Similar to Lemma \ref{Lem L empty}, we will obtain a contradiction. Hence, $|\beta(Q_{1})|\leq9(M+2t)|Q_{1}|$. For any $Q_{i}$ with $2\leq i\leq t$, we can use a similar discussion.
This finishes the proof.
 \hfill$\Box$

\medskip

The general idea of Lemma \ref{symmetry extended} should be attributed to Simonovits \cite{S1}.

\begin{lem}\label{symmetry extended}
 Let $1\leq i_{0}\leq r$ be fixed. Assume that $u_{1},u_{2},...,u_{\tau}$ are symmetric vertices of $G[V_{i_{0}}\cup W]$, and $\left\{u_{1},u_{2},...,u_{\tau}\right\}\subseteq\overline{V}_{i_{0}}$. Let $B$ be the set of vertices in $V(G)-(V_{i_{0}}\cup W)$ which are not adjacent to any one of  $u_{1},u_{2},...,u_{\tau}$. If $u_{1},u_{2},...,u_{\tau}$ are not symmetric in $G$, then there exists one vertex $v\notin B\cup(V_{i_{0}}\cup W)$ which is not adjacent to at least $\frac{1}{2t}\tau$ vertices among  $u_{1},u_{2},...,u_{\tau}$.
\end{lem}

\f{\bf Proof:}  We can assume that $\frac{1}{2t}\tau\geq1$. Otherwise, there is nothing to prove. Then $\tau\geq2t$. For convenience, assume $i_{0}=1$. Since $u_{1},u_{2},...,u_{\tau}$ are symmetric in $G[V_{1}\cup W]$, and $wu_{i}$ is not an edge for any $w\in B$ and $1\leq i\leq\tau$, we see that $u_{1},u_{2},...,u_{\tau}$ are symmetric in $G[B\cup(V_{1}\cup W)]$, too.

Let $Z=V(G)-(B\cup(V_{1}\cup W))$. If $u_{i}$ is adjacent to all the vertices in $Z$ for any $1\leq i\leq\tau$, then $u_{1},u_{2},...,u_{\tau}$ are symmetric in $G$, too. This contradicts the assumption. Thus, without loss of generality, we can assume that $u_{1}$ is not adjacent to $w_{1}\in Z$. Let $G'$ be the graph obtained from $G$ by adding the edge $u_{1}w_{1}$. Since $G\in{\rm SPEX}(n,\mathcal{F})$ and $\rho(G')>\rho(G)$ by Lemma \ref{subgraph}, we see that $G'$ contains a subgraph $F\in \mathcal{F}$. Clearly, $F$ contains the edge $u_{1}w_{1}$.

Let $v_{1},...,v_{j}$ be all the vertices in $Z$ which are  contained in $F$, where $1\leq j\leq t$. Now we show that there are at most $t$ vertices among $u_{1},u_{2},...,u_{\tau}$, which are adjacent to all the vertices $v_{1},...,v_{j}$.
Suppose not.  We can find a vertex $u_{\ell}\notin V(F)$ with $1\leq \ell\leq \tau$, such that $u_{\ell}$ is  adjacent to all the vertices $v_{1},...,v_{j}$. Since $u_{1},u_{2},...,u_{\tau}$ are symmetric in $G[B\cup(V_{1}\cup W)]$, and $u_{\ell}$ is adjacent to all the vertices  $v_{1},...,v_{j}$, we see that all the neighbors of $u_{1}$ in $F$ are adjacent to $u_{\ell}$ in $G$. Let $F'$ be the subgraph of $G$ induced by $(V(F)-\left\{u_{1}\right\})\cup\left\{u_{\ell}\right\}$. Clearly, $F\subseteq F'$. Thus $G$ contains a copy of $F$, a contradiction. Therefore, there are at most $t$ vertices among $u_{1},u_{2},...,u_{\tau}$, which are adjacent to all the vertices $v_{1},...,v_{j}$. Hence, there are at least $\tau-t$ vertices among $u_{1},u_{2},...,u_{\tau}$, which are not adjacent to all the vertices $v_{1},...,v_{j}$. So, we can find a vertex among $v_{1},...,v_{j}$, say $v$, such that $v$ is not adjacent to at least $\frac{\tau-t}{j}\geq\frac{\tau-t}{t}\geq\frac{1}{2t}\tau$ vertices among $u_{1},u_{2},...,u_{\tau}$.
This finishes the proof.
 \hfill$\Box$

\begin{lem}\label{symmetry to G}
 Let $1\leq i_{0}\leq r$ be fixed. Assume that $u_{1},u_{2},...,u_{\tau}$ are symmetric vertices of  $G[\overline{V}_{i_{0}}]$. Then there are at least $(\frac{1}{2t})^{9(M+2t)}2^{-M}\tau$ vertices among $u_{1},u_{2},...,u_{\tau}$ are symmetric in $G$, too.
\end{lem}

\f{\bf Proof:} We can assume that $(\frac{1}{2t})^{9(M+2t)}2^{-M}\tau\geq1$. Otherwise, there is nothing to prove. Then $\tau\geq(2t)^{9(M+2t)}2^{M}$. For convenience, assume $i_{0}=1$. By Lemma \ref{symmetry edge bound}, $|\beta(\left\{u_{i}\right\})|\leq9(M+2t)$ for any $1\leq i\leq\tau$.  For each $u_{i}$, there are at most $2^{M}$  ways of connecting $u_{i}$ to $W$. Hence, there are at least $\tau2^{-M}$ vertices among $u_{1},u_{2},...,u_{\tau}$ are symmetric in $G[V_{1}\cup W]$ too, i.e., connecting to $W$ in the same way. Without loss of generality, assume that $u_{1},u_{2},...,u_{\tau_{0}}$ are the  symmetric vertices in $G[V_{1}\cup W]$, where $\tau_{0}\geq\tau2^{-M}$.

Let $B_{0}$ be the set of vertices in $V(G)-(V_{1}\cup W)$ which are not adjacent to any one of  $u_{1},u_{2},...,u_{\tau_{0}}$. If $u_{1},u_{2},...,u_{\tau_{0}}$ are symmetric in $G$, then the lemma holds as $\tau_{0}\geq\tau2^{-M}$. So, we can assume that $u_{1},u_{2},...,u_{\tau_{0}}$ are not symmetric in $G$. By Lemma \ref{symmetry extended}, there exists one vertex $v_{1}\notin B_{0}\cup(V_{1}\cup W)$, such that $v_{1}$ is not adjacent to at least $\frac{1}{2t}\tau_{0}$ vertices among $u_{1},u_{2},...,u_{\tau_{0}}$. Without loss of generality, let $u_{1},u_{2},...,u_{\tau_{1}}$ be the vertices not adjacent to $v_{1}$, where $\tau_{1}\geq\frac{1}{2t}\tau_{0}\geq\frac{1}{2t}2^{-M}\tau$.

Let $B_{1}$ be the set of vertices in $V(G)-(V_{1}\cup W)$ which are not adjacent to any one of  $u_{1},u_{2},...,u_{\tau_{1}}$. Then $v_{1}\in B_{1}$ and thus $|B_{0}|<|B_{1}|$. If $u_{1},u_{2},...,u_{\tau_{1}}$ are symmetric in $G$, then the lemma holds as $\tau_{1}\geq\frac{1}{2t}2^{-M}\tau$. So, we can assume that $u_{1},u_{2},...,u_{\tau_{1}}$ are not symmetric in $G$. By Lemma \ref{symmetry extended}, there exists one vertex $v_{2}\notin B_{1}\cup(V_{1}\cup W)$, such that $v_{2}$ is not adjacent to at least $\frac{1}{2t}\tau_{1}$ vertices among $u_{1},u_{2},...,u_{\tau_{1}}$. Without loss of generality, let $u_{1},u_{2},...,u_{\tau_{2}}$ be the vertices not adjacent to $v_{2}$, where $\tau_{2}\geq\frac{1}{2t}\tau_{1}\geq(\frac{1}{2t})^{2}2^{-M}\tau$.

Repeat the above process. We can obtain subsets $B_{0}\subset B_{1}\subset\cdots\subset B_{\ell}$, and vertices $u_{1},u_{2},...,u_{\tau_{\ell}}$ which are symmetric in $G[B_{\ell}\cup(V_{1}\cup W)]$, where $\tau_{\ell}\geq(\frac{1}{2t})^{\ell}2^{-M}\tau$ for any $\ell\geq1$. Moreover,  all the vertices in $B_{\ell}$ are not adjacent to any one of $u_{1},u_{2},...,u_{\tau_{\ell}}$. Since $|\beta(\left\{u_{i}\right\})|\leq9(M+2t)$ for any $1\leq i\leq\tau$ and $|B_{\ell}|\geq\ell$, we see that $\ell\leq9(M+2t)$. Since the process must stop at step $\ell\leq9(M+2t)$, we see that $u_{1},u_{2},...,u_{\tau_{\ell}}$ are symmetric in $G$, too. The proof is completed as $\tau_{\ell}\geq(\frac{1}{2t})^{\ell}2^{-M}\tau\geq(\frac{1}{2t})^{9(M+2t)}2^{-M}\tau$.
 \hfill$\Box$

\medskip

\f{\bf Proof of Theorem \ref{Pl-free}.} By Lemma \ref{symmetric}, if a $\left\{P_{t}\right\}$-free graph $H$ has order at least $C(t,k)$, then there are at least $k$ symmetric subgraphs of order at most $D(t,k)$ in $H$. Let $C=C(t,\lceil(2t)^{9(M+2t)}2^{M}t+t\rceil)$ and $D=D(t,\lceil(2t)^{9(M+2t)}2^{M}t+t\rceil)$. Set $M_{0}=\max\left\{(2t)^{9(M+2t)}2^{M+M_{1}}t+M_{1},C\right\}$, where $M_{1}=9(M+2t)Dt+2(t-1)$.

\f{\bf Case 1.} Among $\overline{V}_{1},\overline{V}_{2},...,\overline{V}_{r},$ there exits one, say $\overline{V}_{1}$, such that $G[\overline{V}_{1}]$ contains $t$ non-trivial symmetric subgraphs of order at most $D$.

Let $Q_{1},Q_{2},...,Q_{t}$ be non-trivial symmetric subgraphs of $G[\overline{V}_{1}]$, where $2\leq|Q_{i}|\leq D$ for $1\leq i\leq t$. By Lemma \ref{symmetry edge bound}, $|\beta(Q_{i})|\leq9(M+2t)|Q_{i}|\leq9(M+2t)D$ for any $1\leq i\leq t$. Let $\overline{N}$ be the set of vertices in $\cup_{2\leq i\leq r}\overline{V}_{i}$ which are not adjacent to at least one vertex in $\cup_{1\leq i\leq t}V(Q_{i})$. Then $|\overline{N}|\leq\sum_{1\leq i\leq t}|\beta(Q_{i})|\leq9(M+2t)Dt$. For any $2\leq i\leq r$, let $\overline{V}'_{i}=\overline{V}_{i}-\overline{N}$. Then all the vertices in $\cup_{1\leq i\leq t}V(Q_{i})$ are adjacent to all the vertices in $\cup_{2\leq i\leq r}\overline{V}'_{i}$.

Suppose that $G[\overline{V}'_{i_{0}}]$ has a matching $u_{1}v_{1},u_{2}v_{2},...,u_{t}v_{t}$ for some $2\leq i_{0}\leq r$. Since $Q_{1},Q_{2},...,Q_{t}$ are non-trivial, let $f_{i}g_{i}\in E(Q_{i})$ for any $1\leq i\leq t$. Note that any one of the vertices $u_{1},v_{1},u_{2},v_{2},...,u_{t},v_{t}$ is adjacent to all the vertices $f_{1},g_{1},f_{2},g_{2},...,f_{t},g_{t}$. Similar to Lemma \ref{symmetry edge bound}, we can obtain $tK_{2}\prod tK_{2}\prod T(t(r-2),r-2)\subseteq G$, a contradiction. Hence $G[\overline{V}'_{i}]$ has matching number at most $t-1$ for any $2\leq i\leq r$. Let $2\leq i\leq r$ be any fixed. Let $\overline{V}''_{i}$ be the subset of $\overline{V}'_{i}$ formed by deleting the vertices in a maximum matching in $G[\overline{V}'_{i}]$.
 Then $\overline{V}''_{i}$ is an independent set of $G$, and $|\overline{V}_{i}-\overline{V}''_{i}|\leq|\overline{N}|+2(t-1)\leq9(M+2t)Dt+2(t-1)$. Recall that $M_{1}=9(M+2t)Dt+2(t-1)$. Then $|\overline{V}_{i}-\overline{V}''_{i}|\leq M_{1}$.

 Let $2\leq i\leq r$ be any fixed. Let $w_{1},w_{2},...,w_{\tau}$ be the symmetric vertices in $G$ with $\tau$ maximum, where $w_{j}\in\overline{V}''_{i}$ for any $1\leq j\leq\tau$. We first show that $\tau\geq t$. Considering the ways of connecting the vertices in $\overline{V}''_{i}$ to  $\overline{V}_{i}-\overline{V}''_{i}$, we see that (for large $n$) there are at least $2^{-M_{1}}|\overline{V}''_{i}|\geq (2t)^{9(M+2t)}2^{M}t$  vertices in $\overline{V}''_{i}$ connecting to $\overline{V}_{i}-\overline{V}''_{i}$ in the same way. Thus, these $(2t)^{9(M+2t)}2^{M}t$ vertices are symmetric in $G[\overline{V}_{i}]$. By Lemma \ref{symmetry to G}, there are $t$ vertices (among these $(2t)^{9(M+2t)}2^{M}t$ vertices)  which are symmetric in $G$, too. Therefore, $\tau\geq t$.

Now we show that there are at most $(2t)^{9(M+2t)}2^{M+M_{1}}t-1$ vertices in $\overline{V}''_{i}-\left\{w_{1},w_{2},...,w_{\tau}\right\}$. If not, consider the ways of connecting $(2t)^{9(M+2t)}2^{M+M_{1}}t$ vertices in $\overline{V}''_{i}-\left\{w_{1},w_{2},...,w_{\tau}\right\}$ to  $\overline{V}_{i}-\overline{V}''_{i}$. Recalling that $|\overline{V}_{i}-\overline{V}''_{i}|\leq M_{1}$, there are at least $(2t)^{9(M+2t)}2^{M}t$  vertices in $\overline{V}''_{i}-\left\{w_{1},w_{2},...,w_{\tau}\right\}$ connecting to $\overline{V}_{i}-\overline{V}''_{i}$ in the same way. Thus, these $(2t)^{9(M+2t)}2^{M}t$ vertices are symmetric in $G[\overline{V}_{i}]$. By Lemma \ref{symmetry to G}, there are $t$ vertices (among these $(2t)^{9(M+2t)}2^{M}t$ vertices), say $z_{1},z_{2},...,z_{t}$,  which are symmetric in $G$, too. By the choice of $w_{1},w_{2},...,w_{\tau}$, we see that $w_{1}$ and $z_{1}$ are not symmetric in $G$, i.e., $N_{G}(w_{1})\neq N_{G}(z_{1})$. Without loss of generality, assume $\sum_{v\in N_{G}(w_{1})}x_{v}\geq\sum_{v\in N_{G}(z_{1})}x_{v}$. Let $G'$ be the graph obtained from $G$ by deleting the edges incident with $z_{1}$, and adding all the edges between $z_{1}$ and $N_{G}(w_{1})$. Since $w_{1},w_{2},...,w_{\tau}$ and $z_{1}$ are symmetric in $G'$, we see that if $G'$ contains a $F\subseteq\mathcal{F}$, then $G$ contains a copy of $F$, too. Thus, $G'$ is $\mathcal{F}$-free. But $\rho(G')>\rho(G)$ by Lemma \ref{eigenvector trans}, a contradiction.
Therefore, there are at most $(2t)^{9(M+2t)}2^{M+M_{1}}t-1$ vertices in $\overline{V}''_{i}-\left\{w_{1},w_{2},...,w_{\tau}\right\}$. Let $\overline{V}'''_{i}=\left\{w_{1},w_{2},...,w_{\tau}\right\}$. Then $|\overline{V}''_{i}-\overline{V}'''_{i}|\leq(2t)^{9(M+2t)}2^{M+M_{1}}t$.
Recall that $|\overline{V}_{i}-\overline{V}''_{i}|\leq M_{1}$.  Then $|\overline{V}_{i}-\overline{V}'''_{i}|\leq (2t)^{9(M+2t)}2^{M+M_{1}}t+M_{1}$.

From the above discussion, for any $2\leq i\leq r$, we can obtain a subset $\overline{V}'''_{i}\subseteq\overline{V}_{i}$, such that $|\overline{V}_{i}-\overline{V}'''_{i}|\leq (2t)^{9(M+2t)}2^{M+M_{1}}t+M_{1}$,  and the vertices in $\overline{V}'''_{i}$ are symmetric in $G$. As is shown in Lemma \ref{Main W bound}, $d_{\overline{V}_{j}}(v)\geq\frac{n}{2r}-\frac{r\theta n}{2}$ for any $v\in V_{i}$ and $1\leq i\neq j\leq r$. Moreover, for any $1\leq i\leq r$, $d_{V_{i}}(v)< \theta n$ for any $v\in \overline{V}_{i}$, and $d_{V_{i}}(v)\geq \theta n$ for any $v\in W_{i}$. It is easy to see that $N_{G}(u)=V(G)-\overline{V}_{i}$ for any $u\in \overline{V}'''_{i}$. Now move all the vertices in $\cup_{2\leq i\leq r}W_{i}$ to $V_{1}$. That is to say, let $S_{1}=V_{1}\cup W$, and let $S_{i}= \overline{V}_{i}$ and $S'_{i}=\overline{V}'''_{i}$ for any $2\leq i\leq r$. Then $|S_{i}-S'_{i}|\leq(2t)^{9(M+2t)}2^{M+M_{1}}t+M_{1}\leq M_{0}$ for any $2\leq i\leq r$. Since $\theta>0$ is fixed, $M_{0}$ is a constant (with respective to $n$). Now $V(G)=\cup_{1\leq i\leq r}S_{i}$ is the desired partition.

\f{\bf Case 2.} For any $1\leq i\leq r$, if $Q_{1},Q_{2},...,Q_{t}$ are symmetric subgraphs of $G[\overline{V}_{i}]$ with $|Q_{1}|\leq D$, then these $t$ symmetric subgraphs must be trivial (i.e., $|Q_{1}|=1$).

Let $2\leq i\leq r$ be any fixed. Let $w_{1},w_{2},...,w_{\tau}$ be the symmetric vertices in $G$ with $\tau$ maximum, where $w_{j}\in\overline{V}_{i}$ for any $1\leq j\leq\tau$. We first show that $\tau\geq t$. When $n$ is large, we have $|\overline{V}_{i}|\geq C$. By Lemma \ref{symmetric}, there are at least  $(2t)^{9(M+2t)}2^{M}t+t$ symmetric subgraphs of order at most $D$ in $G[\overline{V}_{i}]$. By the assumption in Case 2, these symmetric subgraphs are all trivial.  By Lemma \ref{symmetry to G}, there are at least $t$ ones among these $(2t)^{9(M+2t)}2^{M}t+t$ vertices are symmetric in $G$, too. Therefore, $\tau\geq t$.
  Let $O$ be the set of vertices in $\overline{V}_{i}$, which are adjacent to $w_{1},w_{2},...,w_{\tau}$. (Since $w_{1},w_{2},...,w_{\tau}$ are symmetric in $G$, another vertex $v$ is adjacent to one of them if and only if $v$ is adjacent to all of them.) Clearly, $|O|\leq t$. Otherwise, there is a path $P_{t}$ inside $\overline{V}_{i}$. Similar to Lemma \ref{symmetry edge bound}, $G$ will contain a copy of $P_{t}\prod T(t(r-1),r-1)$, a contradiction.

Recall that $2\leq i\leq r$ is any fixed. Now we show that there are at most $C-1$ vertices in $\overline{V}_{i}-\left\{w_{1},w_{2},...,w_{\tau}\right\}$. If not, by the definition of $C$, there are at least $(2t)^{9(M+2t)}2^{M}t+t$ symmetric subgraphs in $\overline{V}_{i}-\left\{w_{1},w_{2},...,w_{\tau}\right\}$. Since $|O|\leq t$,  there are at least $(2t)^{9(M+2t)}2^{M}t$ symmetric subgraphs in $\overline{V}_{i}$ among those $(2t)^{9(M+2t)}2^{M}t+t$ subgraphs. By the assumption of Case 2, these $(2t)^{9(M+2t)}2^{M}t$ symmetric subgraphs  in $\overline{V}_{i}$  are trivial.  By Lemma \ref{symmetry to G}, there are at least $t$ ones among these $(2t)^{9(M+2t)}2^{M}t$ symmetric vertices  in $\overline{V}_{i}$ are symmetric in $G$, too. Let $z_{1},z_{2},...,z_{t}$ be the $t$ symmetric vertices  in $G$. By the choice of $w_{1},w_{2},...,w_{\tau}$,  we will obtain a contradiction similar to Case 1. Hence there are at most $C-1$ vertices in $\overline{V}_{i}-\left\{w_{1},w_{2},...,w_{\tau}\right\}$. Let $\overline{V}'_{i}=\left\{w_{1},w_{2},...,w_{\tau}\right\}$. Then $|\overline{V}_{i}-\overline{V}'_{i}|\leq C\leq M_{0}$. let $S_{1}=V_{1}\cup W$, and let $S_{i}= \overline{V}_{i}$ and $S'_{i}=\overline{V}'_{i}$ for any $2\leq i\leq r$. Then $|S_{i}-S'_{i}|\leq M_{0}$. The rest of the proof is very similar to Case 1, and thus omitted. This completes the proof. \hfill$\Box$

\medskip

\f{\bf Proof of Theorem \ref{matching-free}.} Recall that, if $F_{0}\subseteq tK_{2}\prod T(t(r-1),r-1)$ for some $F_{0}\in\mathcal{F}$, then $\mathcal{F}$ satisfies property $(P)$. We can follow the proof of Theorem \ref{Pl-free}. As in Theorem \ref{Pl-free}, we see that any two vertices from distinct parts of $S_{1},S_{2},...,S_{r}$ are adjacent in $G$. Thus, the matching number of $G[S_{1}]$ is  at most $t$. Let $R$ be the set of vertices in $S_{1}$ which are the end vertices of a maximum matching of $G[S_{1}]$. Then $|R|\leq 2t$, and $S_{1}-R$ is an independent set in $G$. (Clearly, $W\subseteq R$ as $W\subseteq S_{1}$.) Let $\overline{S}_{1}=S_{1}-R$.

Let $w_{1},w_{2},...,w_{\tau}$ be the symmetric vertices in $G$ with $\tau$ maximum, where $w_{j}\in\overline{S}_{1}$ for any $1\leq j\leq\tau$. We first show that $\tau\geq t$. For large $n$, we have $2^{-2t}|\overline{S}_{1}|\geq(2t)^{9(M+2t)}2^{M}t$. Thus, there are at least  $(2t)^{9(M+2t)}2^{M}t$  vertices in $\overline{S}_{1}$ connecting to $R$ in the same way. That is to say, these $(2t)^{9(M+2t)}2^{M}t$ vertices are symmetric in $G[S_{1}]$. By Lemma \ref{symmetry to G}, there are $t$ vertices among them, which are symmetric in $G$, too. Therefore, $\tau\geq t$.

In the same way as the proof of Theorem \ref{Pl-free},  we can show that there are at most $(2t)^{9(M+2t)}2^{M+2t}t-1$ vertices in $\overline{S}_{1}-\left\{w_{1},w_{2},...,w_{\tau}\right\}$. Let $S'_{1}=\left\{w_{1},w_{2},...,w_{\tau}\right\}$. Then $|S_{1}-S'_{1}|\leq(2t)^{9(M+2t)}2^{M+2t}t+2t\leq M_{0}$. We also remove the vertices in $W$ out of $S_{1}$. By the definition of $W$, we see that $d_{S_{i}}(v)\geq \frac{1}{2}\theta n$ for any $v\in W$ and $1\leq i\leq r$.  Hence, each vertex $v\in W$ is adjacent to all the vertices in $\cup_{1\leq i\leq r}\overline{S}_{i}$. Clearly, $|W|\leq q-1$, since $F\subseteq \overline{K_{q}}\prod T(tr,r)$ for some $F\subseteq\mathcal{F}$.

From the above discussion, we see that $V(G)=W\cup(\cup_{1\leq i\leq r}S_{i})$, where $|W|\leq q-1$. Moreover, for any $1\leq i\leq r$, there is a subset $S'_{i}\subseteq S_{i}$ (with $|S_{i}-S'_{i}|\leq M_{0}$) such that $N_{G}(u)=V(G)-S_{i}$ for any $u\in S'_{i}$. Since $\theta>0$ is fixed, $M_{0}$ is a constant. Thus, we can choose a small constant $\sigma>0$, such that $(\frac{r-1}{r}-6\sigma)(M_{0}+1)>(\frac{r-1}{r}+5r\sigma)M_{0}$ (i.e., $\sigma<\frac{r-1}{r(5r+6)M_{0}}$). By Lemma \ref{mini degree}, we can assume that $\delta(G)>(1-\frac{1}{r}-\sigma)n$ for sufficiently large $n$. Since for any $1\leq i\leq r$, each vertex in $S'_{i}$ has degree $n-|S_{i}|$, we have $n-|S_{i}|>(1-\frac{1}{r}-\sigma)n$. This implies that $|S_{i}|<(\frac{1}{r}+\sigma)n$ for any $1\leq i\leq r$. Since $\sum_{1\leq i\leq r}|S_{i}|=n-|W|>(1-\sigma)n$ for large $n$. It follows that $|S_{i}|>(1-\sigma)n-(r-1)(\frac{1}{r}+\sigma)n=(\frac{1}{r}-r\sigma)n$. Hence $(\frac{1}{r}-r\sigma)n<|S_{i}|<(\frac{1}{r}+\sigma)n$ for any $1\leq i\leq r$.

Recall that $x_{u^{*}}=\max_{v\in V(G)}x_{v}=1$ and $x_{v^{*}}=\max_{v\in V(G)-W}x_{v}$.

\medskip

\f{\bf Claim 1.} $1-\sigma\leq x_{v}\leq1$ for any $v\in W$, and $\frac{r-1}{r}-6\sigma\leq x_{u}\leq\frac{r-1}{r}+5r\sigma$ for any $u\in V(G)-W=\cup_{1\leq i\leq r}S_{i}$.

\medskip

\f{\bf Proof of Claim 1.} Recall that $\rho(G)\geq\frac{r-1}{r}n-\frac{r}{4n}$. For any $v\in W$, since
\begin{equation}
\begin{aligned}
\rho(G)x_{v}-\rho(G)x_{u^{*}}&=(\sum_{u\in N_{G}(v)}x_{u})-(\sum_{u\in N_{G}(u^{*})}x_{u})\\
&\geq-(\sum_{u\in W}x_{u})-(\sum_{u\in \cup_{1\leq i\leq r}(S_{i}-S'_{i})}x_{u})\\
&\geq-(q-1)-M_{0}r,
\end{aligned}\notag
\end{equation}
we have that $x_{v}\geq1-\frac{q-1+M_{0}r}{\rho(G)}\geq1-\sigma$. Hence $1-\sigma\leq x_{v}\leq1$ for any $v\in W$.

Using a quite similar way to Lemma \ref{entry lower}. we can show that $x_{u}\geq(1-5\sigma)x_{v^{*}}$ and $x_{u}\geq\frac{r-1}{r}-6\sigma$ for any $v\in V(G)-W$. For the upper bound, it suffices to show $x_{v^{*}}\leq\frac{r-1}{r}+5r\sigma$.
 Assume that $v^{*}\in S_{i_{0}}$. Note that $|S'_{i_{0}}|\geq(\frac{1}{r}-r\sigma)n$.
Since
\begin{equation}
\begin{aligned}
\rho(G)x_{v^{*}}-\rho(G)x_{u^{*}}&=(\sum_{v\in N_{G}(v^{*})}x_{v})-(\sum_{v\in N_{G}(u^{*})}x_{v})\\
&\leq(\sum_{v\in W}x_{v})+(\sum_{v\in \cup_{1\leq i\leq r}(S_{i}-S'_{i})}x_{v})-(\sum_{v\in S'_{i_{0}}}x_{v})\\
&\leq(q-1)+M_{0}r-(1-5\sigma)x_{v^{*}}|S'_{i_{0}}|\\
&\leq(q-1)+M_{0}r-(1-5\sigma)(\frac{1}{r}-r\sigma)nx_{v^{*}}\\
&\leq(q-1)+M_{0}r-(\frac{1}{r}-3r\sigma)nx_{v^{*}},
\end{aligned}\notag
\end{equation}
we have that $$x_{v^{*}}\leq1+\frac{q-1+M_{0}r}{\rho(G)}-\frac{(\frac{1}{r}-3r\sigma)n}{\rho(G)}x_{v^{*}}\leq1+\sigma-\frac{1-4r^{2}\sigma}{r-1}x_{v^{*}}.$$
It follows that $x_{v^{*}}\leq\frac{1+\sigma}{\frac{r-4r^{2}\sigma}{r-1}}\leq\frac{r-1}{r}+5r\sigma$, as desired. \hfill$\Box$

Now we show that $|W|=q-1$. If $q=1$, there is nothing to prove. Let $q\geq2$. Suppose $|W|<q-1$. Let $u$ be a vertex in $S'_{1}$. Let $G'$ be the graph obtained from $G$ by deleting all the edges in $E(W)\cup(\cup_{1\leq i\leq r}E(S_{i}-S'_{i}))$, and adding the edges between $u$ and $S_{1}-\left\{u\right\}$. Clearly, $G'$ is $\mathcal{F}$-free, since $G'$ is a subgraph of $\overline{K_{q-1}}\prod T(rn,r)$.
However, using Claim 1,
\begin{equation}
\begin{aligned}
\mathbf{x}^{T}(\rho(G')-\rho(G))\mathbf{x}&\geq\mathbf{x}^{T}(A(G')-A(G))\mathbf{x}\\
&\geq2(\sum_{v\in S_{1}-\left\{u\right\}}x_{u}x_{v})-2(\sum_{wv\in E(W)\cup(\cup_{1\leq i\leq r}E(S_{i}-S'_{i}))}x_{w}x_{v})\\
&\geq2(\frac{r-1}{r}-6\sigma)^{2}(|S_{1}|-1)-2|E(W)\cup(\cup_{1\leq i\leq r}E(S_{i}-S'_{i}))|\\
&\geq2(\frac{r-1}{r}-6\sigma)^{2}(|S_{1}|-1)-2((q-1)^{2}+rM_{0}^{2})\\
&\geq2(\frac{r-1}{r}-6\sigma)^{2}(\frac{1}{r}-2r\sigma)n-2(q-1)^{2}-2rM_{0}^{2}\\
&>0.
\end{aligned}\notag
\end{equation}
It follows  that $\rho(G')>\rho(G)$. This contradicts the fact that $G\in {\rm SPEX}(n,\mathcal{F})$. Therefore, $|W|=q-1$.

Set $n_{i}=|S_{i}|$ for any $1\leq i\leq r$.
Without loss of generality, assume $n_{1}\geq n_{2}\geq\cdots\geq n_{r}$.
It remains to show that $n_{i}=\lfloor\frac{n-q+1}{r}\rfloor$ or $\lceil\frac{n-q+1}{r}\rceil$ for each $1\leq i\leq r$.
Suppose not. Then $n_{1}\geq n_{r}+2$. Recall that $|S_{i}-S'_{i}|\leq M_{0}$ for any $1\leq i\leq r$, where $M_{0}$ can be regarded as an integer.
Let $A$ be an independent subset of $S_{1}$ with $|A|=n_{1}-(M_{0}+1)$, and let $B$ be an independent subset of $S_{r}$ with $|B|=n_{r}-M_{0}$. Since $n_{1}\geq n_{2}+2$, we have $|A|>|B|>0$.
Since $\frac{r-1}{r}-6\sigma\leq x_{u}\leq\frac{r-1}{r}+5r\sigma$ for any $u\in \cup_{1\leq i\leq r}S_{i}$, we have that
 $$\sum_{v\in S_{1}-A}x_{v}\geq(\frac{r-1}{r}-6\sigma)(M_{0}+1)>(\frac{r-1}{r}+5r\sigma)M_{0}\geq\sum_{v\in S_{r}-B}x_{v}.$$
Let $U=V(G)-(A\cup B)$. Let $U_{1}=(S_{r}-B)\cup (V(G)-(S_{1}\cup S_{r}))$ and $U_{2}=(S_{1}-A)\cup (V(G)-(S_{1}\cup S_{r}))$. Clearly, $N_{U}(v)=U_{1}$ for any $v\in A$, and $N_{U}(v)=U_{2}$ for any $v\in B$. Moreover, $\sum_{v\in U_{1}}x_{v}<\sum_{v\in U_{2}}x_{v}$, since $\sum_{v\in S_{1}-A}x_{v}>\sum_{v\in S_{r}-B}x_{v}$ as above.
Let $G''$ be the graph obtained from $G$ by deleting all the edges between $A$ and $U_{1}$, and between $B$ and $U_{2}$, and adding all the edges between $A$ and $U_{2}$, and between $B$ and $U_{1}$. By Lemma \ref{eigenvector trans var}, we have $\lambda_{1}(G'')>\lambda_{1}(G)$. Clearly, $G''$ is $\mathcal{F}$-free, since $G$ is $\mathcal{F}$-free. But, this contradicts the fact that $G\in {\rm SPEX}(n,\mathcal{F})$. Hence $n_{1}\leq n_{r}+1$. It follows that $n_{i}=\lfloor\frac{n-q+1}{r}\rfloor$ or $\lceil\frac{n-q+1}{r}\rceil$ for each $1\leq i\leq r$.  This completes the proof. \hfill$\Box$

\section{Proofs of Theorem \ref{application 1} and Theorem \ref{application 2}}

\f{\bf Proof of Theorem \ref{application 1}.}
Recall that $1\leq s\leq\frac{m}{2}$, $r\geq2$, and $F=(sK_{2}\cup\overline{K_{m-2s}})\prod T(m(r-1),r-1)$. Clearly, $\chi(F)=r+1$. Let $t=|F|$. It is easy to see that $F\subseteq \overline{K_{\ell}}\prod T(tr,r)$ for any $\ell\geq s$, and $F\nsubseteq K_{\ell}\prod T(tr,r)$ for any $\ell\leq s-1$. This implies that $q(F)=s$. For large $n$, let $G\in{\rm SPEX}(n,\left\{F\right\})$. Since $F\subseteq(t K_{2}\prod T(t(r-1),r-1)$, by Theorem \ref{matching-free}, we see that $G$ has a partition $V(G)=W\cup(\cup_{1\leq i\leq r}S_{i})$ satisfying the following property:\\
$(i)$ $|W|=s-1$ and $|S_{i}|=\lfloor\frac{n-s+1}{r}\rfloor$ or $\lfloor\frac{n-s+1}{r}\rfloor$ for any $1\leq i\leq r$;\\
$(ii)$ for any $1\leq i\leq r$, there exists $S'_{i}\subseteq S_{i}$ such that $|S_{i}-S'_{i}|=O(1)$, and $N_{G}(v)=V(G)-S_{i}$ for any $v\in S'_{i}$.

Now we show that there are no edges inside $S_{i}$ for any $1\leq i\leq r$. If not, assume that $uv$ is an edge in $S_{i_{0}}$, where $1\leq i_{0}\leq r$. Clearly, $G[W\cup S_{i_{0}}]$ contains a matching of $s$ edges, say $uv,u_{1}v_{1},...,u_{s-1}v_{s-1}$. Let $w_{1},w_{2},...,w_{m-2s}$ be another $m-2s$ vertices in $S_{i_{0}}$. Note that any two vertices from distinct parts: $W,S'_{j}$ for all $1\leq j\leq r$, are adjacent in $G$,  and all the vertices in $\cup_{1\leq j\neq i_{0}\leq r}S'_{j}$ are adjacent to all the vertices in $\left\{u,v,u_{1},v_{1},...,u_{s-1},v_{s-1}\right\}\cup\left\{w_{1},w_{2},...,w_{m-2s}\right\}$. Thus, $G$ contains a copy of $(sK_{2}\cup\overline{K_{m-2s}})\prod T(m(r-1),r-1)$, a contradiction. Hence, there are no edges inside $S_{i}$ for any $1\leq i\leq r$.

From the above discussion,  we see that $G$ is a spanning subgraph of $K_{s-1}\prod T(n-s+1,r)$. Since $G\in{\rm SPEX}(n,\left\{F\right\})$, by Lemma \ref{subgraph}, we have $G=K_{s-1}\prod T(n-s+1,r)$. This completes the proof. \hfill$\Box$

\medskip

\f{\bf Proof of Theorem \ref{application 2}.}
Recall that $m=p(k+1)+h$, where $p\geq2$ and $1\leq h\leq k$.  For large $n$, let $G$ be a graph in  ${\rm SPEX}(n,\left\{C^{k}_{m}\right\})$. We can assume that $h=ap+b$, where $a\geq0$ and $1\leq b\leq p$. Then, it is easy to see that $b=h-\lfloor\frac{h-1}{p}\rfloor p$.
Let $v_{1},v_{2},...,v_{m}$ be the vertices of $C_{m}$, where $v_{i}$ is adjacent to $v_{i+1}$ for any $1\leq i\leq m-1$, and $v_{1}$ is adjacent to $v_{m}$. We can assume that $C^{k}_{m}$ is obtained from $C_{m}$ by adding edges $uv$ whenever $d_{C_{m}}(u,v)\leq k$. Let $\sigma$ be the independence number of $C^{k}_{m}$ (i.e., the maximum cardinality of an independent set of $C^{k}_{m}$). Since any $k+1$ successive vertices among $\left\{v_{1},v_{2},...,v_{m}\right\}$ induce a clique in $C^{k}_{m}$, we have
$\sigma=\lfloor\frac{m}{k+1}\rfloor=p$. Set $r=k+\lceil\frac{h}{p}\rceil$. Then $m=pr+b$.

Now we show that $\chi(C^{k}_{m}-U)\geq r+1$ for any subset $U\subseteq V(C^{k}_{m})$ with $|U|=b-1$. Since an independent set in  $C^{k}_{m}-U$ is also an independent set in  $C^{k}_{m}$, from $\sigma=p$ we see that $$\chi(C^{k}_{m}-U)\geq\lceil\frac{m-b+1}{p}\rceil=k+1+\lceil\frac{h}{p}\rceil=r+1,$$
as desired.

Now we show that $C^{k}_{m}$ is a subgraph of $F=(bK_{2}\cup\overline{K_{m-2b}})\prod T(m(r-1),r-1)$. We first consider the case $1\leq b\leq p-1$.
For any $1\leq i\leq r$, let
$$S_{i}=\left\{v_{i},v_{i+(r+1)},v_{i+2(r+1)},...,v_{i+b(r+1)},v_{i+b(r+1)+r},v_{i+b(r+1)+2r},...,v_{i+b(r+1)+(p-b-1)r}\right\},$$
Let $W=V(G)-\cup_{1\leq i\leq r}S_{i}$.
 Clearly, $S_{i}$ is an independent set for any $1\leq i\leq r$, and
 $$W=\left\{v_{r+1},v_{r+1+(r+1)},...,v_{r+1+(b-1)(r+1)}\right\}.$$
Moreover, $G[W\cup S_{r}]$ consists of a matching with $b$ edges $$v_{r}v_{r+1},v_{r+(r+1)}v_{r+1+(r+1)},...,v_{r+(b-1)(r+1)}v_{r+1+(b-1)(r+1)},$$
 and $p-b$ isolated vertices $$v_{r+1+(b-1)(r+1)+r},v_{r+1+(b-1)(r+1)+2r},...,v_{r+1+(b-1)(r+1)+(p-b)r}.$$
 Thus,  $C^{k}_{m}$ is a subgraph of $F=(bK_{2}\cup\overline{K_{m-2b}})\prod T(m(r-1),r-1)$.

 It remains the case $b=p$. For any $1\leq i\leq r+1$, let
$$S_{i}=\left\{v_{i},v_{i+(r+1)},v_{i+2(r+1)},...,v_{i+(p-1)(r+1)}\right\}.$$
Clearly, $S_{i}$ is an independent set for any $1\leq i\leq r+1$.
Moreover, $G[S_{r}\cup S_{r+1}]$ consists of a matching with $p$ edges $$v_{r}v_{r+1},v_{r+(r+1)}v_{r+1+(r+1)},...,v_{r+(p-1)(r+1)}v_{r+1+(p-1)(r+1)}.$$
 Thus,  $C^{k}_{m}$ is a subgraph of $F=(pK_{2}\cup\overline{K_{m-2p}})\prod T(m(r-1),r-1)$.
Hence, in either case, we see that $C^{k}_{m}$ is a subgraph of $(bK_{2}\cup\overline{K_{m-2b}})\prod T(m(r-1),r-1)$.

By the above discussion, $\chi(C^{k}_{m}-U)\geq r+1$ for any subset $U\subseteq V(C^{k}_{m})$ with $|U|=b-1$, and $C^{k}_{m}\subseteq (bK_{2}\cup\overline{K_{m-2b}})\prod T(m(r-1),r-1)$. It follows that $\chi(C^{k}_{m})=r+1$, and $C^{k}_{m}\not\subseteq K_{b-1}\prod T(n-b+1,r)$.

$(i)$ By a result of Simonovits \cite{S,S1}, for large $n$, $K_{b-1}\prod T(n-b+1,r)$ is the unique graph in
${\rm EX}(n,\left\{(bK_{2}\cup\overline{K_{m-2b}})\prod T(m(r-1),r-1)\right\})$. It follows that $K_{b-1}\prod T(n-b+1,r)$ is the unique graph in ${\rm EX}(n,\left\{C^{k}_{m}\right\})$.

$(ii)$  By Theorem \ref{application 1}, for large $n$, $K_{b-1}\prod T(n-b+1,r)$  is  the unique graph in
${\rm SPEX}(n,\left\{(bK_{2}\cup\overline{K_{m-2b}})\prod T(m(r-1),r-1)\right\})$. It follows that $K_{b-1}\prod T(n-b+1,r)$ is the unique graph in ${\rm SPEX}(n,\left\{C^{k}_{m}\right\})$. This completes the proof. \hfill$\Box$

\medskip

\f{\bf Note.} This is an updated version, since there is a symbol conflict  at Claim 1 of the proof of Theorem \ref{matching-free} in the original manuscript.

\medskip

\f{\bf Declaration of competing interest}

\medskip

There is no conflict of interest.

\medskip

\f{\bf Data availability statement}

\medskip

No data was used for the research described in the article.

\end{document}